\documentclass{amsart}

\usepackage{amsmath}
\usepackage{amsfonts}
\usepackage{amssymb,enumerate}
\usepackage{amsthm}
\usepackage[all]{xy}
\usepackage{hyperref}

\newtheorem{lem}{Lemma}[section]
\newtheorem{cor}[lem]{Corollary}
\newtheorem{prop}[lem]{Proposition}
\newtheorem{thm}[lem]{Theorem}

\newtheorem{Defn}[lem]{Definition}
\newtheorem{Ex}[lem]{Example}
\newtheorem{Question}[lem]{Question}
\newtheorem{Property}[lem]{Property}
\newtheorem{Properties}[lem]{Properties}
\newtheorem{Discussion}[lem]{Remark}
\newtheorem{Construction}[lem]{Construction}
\newtheorem{Notation}[lem]{Notation}
\newtheorem{Fact}[lem]{Fact}
\newtheorem{Notationdefinition}[lem]{Definition/Notation}
\newtheorem{Remarkdefinition}[lem]{Remark/Definition}
\newtheorem{Subprops}{}[lem]
\newtheorem{Para}[lem]{}

\newenvironment{defn}{\begin{Defn}\rm}{\end{Defn}}
\newenvironment{ex}{\begin{Ex}\rm}{\end{Ex}}
\newenvironment{question}{\begin{Question}\rm}{\end{Question}}

\newenvironment{properties}{\begin{Properties}\rm}{\end{Properties}}
\newenvironment{props}{\begin{Properties}\rm}{\end{Properties}}

\newenvironment{fact}{\begin{Fact}\rm}{\end{Fact}}
\newenvironment{notationdefinition}{\begin{Notationdefinition}\rm}{\end{Notationdefinition}}

\newenvironment{subprops}{\begin{Subprops}\rm}{\end{Subprops}}
\newenvironment{para}{\begin{Para}\rm}{\end{Para}}
\newenvironment{disc}{\begin{Discussion}\rm}{\end{Discussion}}

\newcommand{\D}{\mathcal{D}}
\newcommand{\catdb}{\D_{\mathrm{b}}}

\newcommand{\cat}[1]{\mathcal{#1}}

\newcommand{\cata}{\cat{A}}

\newcommand{\catb}{\cat{B}}

\newcommand{\catac}{\cat{A}_C}
\newcommand{\catab}{\cat{A}_B}
\newcommand{\catbc}{\cat{B}_C}

\newcommand{\catbb}{\cat{B}_B}

\newcommand{\pd}{\operatorname{pd}}

\newcommand{\id}{\operatorname{id}}     
\newcommand{\fd}{\operatorname{fd}}

\newcommand{\cidim}{\mathrm{CI}\text{-}\!\dim}

\newcommand{\ciid}{\mathrm{CI}\text{-}\!\id}

\newcommand{\depth}{\operatorname{depth}}       
 
\newcommand{\amp}{\operatorname{amp}}

\newcommand{\cmd}{\operatorname{cmd}}

\newcommand{\ext}{\operatorname{Ext}}   
\newcommand{\rhom}{\mathbf{R}\!\operatorname{Hom}}      
\newcommand{\lotimes}{\otimes^{\mathbf{L}}}
\newcommand{\HH}{\operatorname{H}}
\newcommand{\Hom}{\operatorname{Hom}}   

\newcommand{\spec}{\operatorname{Spec}}
\newcommand{\s}{\ol{\mathfrak{S}}}
\newcommand{\tor}{\operatorname{Tor}}

\newcommand{\shift}{\mathsf{\Sigma}}

\newcommand{\Pic}{\operatorname{Pic}}
\newcommand{\dpic}{\operatorname{DPic}}

\newcommand{\Ker}{\operatorname{Ker}}

\newcommand{\ideal}[1]{\mathfrak{#1}}
\newcommand{\m}{\ideal{m}}

\newcommand{\p}{\ideal{p}}
\newcommand{\q}{\ideal{q}}
\newcommand{\fa}{\ideal{a}}
\newcommand{\fb}{\ideal{b}}

\newcommand{\comp}[1]{\widehat{#1}}
\newcommand{\ol}{\overline}

\newcommand{\ass}{\operatorname{Ass}}
\newcommand{\supp}{\operatorname{Supp}}

\newcommand{\Min}{\operatorname{Min}}

\newcommand{\zz}{\mathbb{Z}}

\newcommand{\from}{\leftarrow}
\newcommand{\xra}{\xrightarrow}
\newcommand{\xla}{\xleftarrow}

\newcommand{\res}{\xra{\simeq}}

\newcommand{\vf}{\varphi}

\newcommand{\y}{\mathbf{y}}

\newcommand{\x}{\mathbf{x}}

\renewcommand{\geq}{\geqslant}
\renewcommand{\leq}{\leqslant}
\renewcommand{\ker}{\Ker}
\renewcommand{\hom}{\Hom}

\begin{document}

\bibliographystyle{amsplain}

\title{Relations between semidualizing complexes}

\author[A.\ J.\ Frankild]{Anders J.~Frankild}
\address{Anders J.~Frankild, University of Copenhagen, Institute for Mathematical 
Sciences, Department of Mathematics, 
Universitetsparken 5, 2100 K\o benhavn, Denmark}
\thanks{This work was completed after the untimely passing of Anders J.\ Frankild
on 10 June 2007.}

\author[S.\ Sather-Wagstaff]{Sean Sather-Wagstaff}
\address{Sean Sather-Wagstaff, Department of Mathematics,
300 Minard Hall,
North Dakota State University,
Fargo, North Dakota 58105-5075, USA}
\email{Sean.Sather-Wagstaff@ndsu.edu}
\urladdr{http://math.ndsu.nodak.edu/faculty/ssatherw/}

\author[A.\ Taylor]{Amelia Taylor}
\address{Amelia Taylor, Colorado College
14 E. Cache La Poudre St.
Colorado Springs, CO 80903, USA}
\email{amelia.taylor@coloradocollege.edu}
\urladdr{http://faculty1.coloradocollege.edu/~ataylor}



\keywords{Auslander classes, Bass classes, complete intersection dimensions, CI-dimensions,
Gorenstein dimensions, G-dimensions, semidualizing complexes, tilting complexes}
\subjclass[2000]{13D05, 13D07, 13D25, 13H10}

\begin{abstract}
We study the following question: Given two semidualizing complexes $B$ and $C$ over a commutative noetherian ring
$R$, does the vanishing of 
$\ext^n_R(B,C)$ for  $n\gg 0$ imply that $B$ is $C$-reflexive?  
This question is
a natural generalization of one studied by Avramov, Buchweitz, and \c{S}ega.
We begin by providing conditions
equivalent to $B$ being $C$-reflexive,
each of which is slightly stronger than the condition
$\ext^n_R(B,C)=0$ for all $n\gg 0$.
We introduce and investigate
an equivalence relation $\approx$ on the set of isomorphism
classes of semidualizing complexes. This relation is defined in terms of a natural
action of the derived Picard group and is well-suited for the study of semidualizing complexes
over nonlocal rings.  We identify numerous alternate characterizations
of this relation, 
each of which includes the condition
$\ext^n_R(B,C)=0$ for all $n\gg 0$.
Finally, we answer our original
question in some special cases.
\end{abstract}

\maketitle


\section{Introduction} \label{sec01}

Given a dualizing complex $D$ for a commutative noetherian ring $R$, 
cohomological properties of $D$ often translate to ring-theoretic properties of  $R$.
For example,  when $R$ is local, if
$\ext^n_R(D,R)=0$ for $n\gg 0$ and the natural 
evaluation morphism
$D\lotimes_R\rhom_R(D,R)\to R$ is an isomorphism in the
derived category $\D(R)$, then $R$ is Gorenstein.
Recently, Avramov, Buchweitz, and \c{S}ega~\cite{avramov:edcrcvct} 
investigated the following potential extensions of this fact.

\begin{para} \label{ABSq}
Let $R$ be a local ring admitting a dualizing complex $D$ such that
$\inf(D)=0$.
\\[1mm]
\noindent\textbf{Question.}  If $\ext^n_R(D,R)=0$ for $(\dim(R)+1)$
consecutive values of  $n\geq 1$, must $R$ be Gorenstein?
\\[1mm]
\noindent\textbf{Conjecture.}  If $\ext^{n}_R(D,R)=0$ for all $n\geq 1$, then $R$ is Gorenstein.
\end{para}

\noindent This paper is concerned with a version of~\eqref{ABSq} for
\emph{semi}dualizing complexes\footnote{While 
this paper is written in the language of complexes and derived categories,
our results specialize readily to the 
case of semidualizing \emph{modules}. 
For a discussion of the translation from complexes
to modules, see~\cite[(4.10)]{christensen:scatac} 
and~\cite[(2.1)]{frankild:rrhffd}.}.

Semidualizing complexes were introduced by 
Avramov and Foxby~\cite{avramov:rhafgd} in a special case
for use in studying local ring homomorphisms,
and by Christensen~\cite{christensen:scatac} in general.
For example, a dualizing complex is semidualizing, as is a free module of rank 1.
Each semidualizing complex $C$ yields a duality theory or,
more specifically, a notion of $C$-reflexivity with properties similar
to those for reflexivity with respect to $D$ or $R$;
see Section~\ref{secintro}  for background material.

Our version of~\eqref{ABSq} for this setting is contained in the
following list of questions.  Specifically,
an affirmative answer to Question~\ref{tach03}\eqref{tach03item1}
would yield an affirmative answer to the question in~\eqref{ABSq},
and an affirmative answer to Question~\ref{tach03}\eqref{tach03item2} would
establish the conjecture in~\eqref{ABSq};
see Remark~\ref{sdc03}. 
Also, note that Example~\ref{local01} shows the need for the local hypothesis 
in Question~\ref{tach03}\eqref{tach03item1}.

\begin{question} \label{tach03}
Let $B$ and $C$ be semidualizing $R$-complexes.
\begin{enumerate}[\quad(a)]
\item \label{tach03item1}
If $\ext^n_R(B,C)=0$ 
for $(\dim(R)+1)$ consecutive values of $n>\sup(C)-\inf(B)$
and if $R$ is local,
must $B$ be $C$-reflexive?
\item \label{tach03item2}
If $\ext^n_R(B,C)=0$ 
for all  $n>\sup(C)-\inf(B)$,
must $B$ be $C$-reflexive?
\item \label{tach03item3}
If $\ext^n_R(B,C)=0$ for $n\gg 0$,
must $B$ be $C$-reflexive?
\end{enumerate}
\end{question}

Our results come in three types.  
The results of 
Section~\ref{sec03} are of the first type:
Assuming $\ext^n_R(B,C)=0$ for $n\gg 0$ \emph{and a bit more},
we show that $B$ is $C$-reflexive. The following is one such
result; its proof is in~\ref{semidual2}.
Note that each of the conditions~\eqref{semidual2item2}--\eqref{semidual2item4}
includes the condition $\ext^n_R(B,C)=0$ for $n\gg 0$.

\begin{thm} \label{intthmA}
Let  $B$ and $C$ be semidualizing $R$-complexes.
The following conditions are equivalent:
\begin{enumerate}[\quad\rm(i)]
\item \label{semidual2item1}
$B$ is $C$-reflexive;
\item \label{semidual2item2}
$\rhom_R(B,C)$ is semidualizing;
\item \label{semidual2item3}
$\rhom_R(B,C)$ is $C$-reflexive;
\item \label{semidual2item4}
$C$ is in the Bass class $\catbb(R)$.
\end{enumerate}
\end{thm}

\noindent
We also prove analogous results for tensor products motivated by the
corresponding version of Question~\ref{tach03} found in
Question~\ref{tpq02}. 

Before discussing our second type of results, 
we describe a new equivalence relation
on the set of isomorphism classes of semidualizing $R$-complexes:
We write $[B]\approx [C]$  if there
is a tilting $R$-complex $P$ such that $B\simeq P\lotimes_RC$.
Here, a \emph{tilting $R$-complex} is a semidualizing $R$-complex of 
finite projective dimension.
When $R$ is local, the only tilting $R$-complexes are
those of the form $\shift^nR$,
so in this case
$[B]\approx [C]$ if and only if
$B$ and $C$ are isomorphic up to shift in $\D(R)$.
Hence, our new relation recovers the more standard relation
over a local local ring while also being particularly well-suited for
the nonlocal setting.
Section~\ref{sec05} contains our treatment of tilting complexes
and the basics of this relation.

Section~\ref{sec04} contains our results of the second type:
Assuming $\ext^n_R(B,C)=0$ for $n\gg 0$ \emph{and a lot more},
we show that $[B]\approx [C]$.
For instance, we prove the next result in~\ref{sdc04}.
As with Theorem~\ref{intthmA}, observe that
each of the conditions~\eqref{sdc04item6}--\eqref{sdc04item5}
includes the condition $\ext^n_R(B,C)=0$ for $n\gg 0$.

\begin{thm} \label{intthmB}
Let $B$ and $C$ be semidualizing $R$-complexes.
The 
following conditions are equivalent:
\begin{enumerate}[\quad\rm(i)]
\item  \label{sdc04item1}
$[B]\approx [C]$;
\item  \label{sdc04item6}
$\rhom_R(B, C)$ is a tilting $R$-complex;
\item  \label{sdc04item2}
$\rhom_R(B, C)$ has finite projective dimension;
\item  \label{sdc04item3}
$\rhom_R(B, C)$ has finite complete intersection dimension;
\item  \label{sdc04item4}
There is an equality of Bass classes $\catbb(R)=\catbc(R)$;
\item  \label{sdc04item5}
$B\in\catbc(R)$ and $C\in\catbb(R)$.
\end{enumerate}
\end{thm}

The last two sections contain our results of the third type: analogues of results
of Avramov, Buchweitz, and \c{S}ega~\cite{avramov:edcrcvct}.
In Section~\ref{sec02} we consider the question of 
when the vanishing assumptions in Question~\ref{tach03}
guarantee that $R$ is Cohen-Macaulay.  This lays some of the foundation
for the results of Section~\ref{sec06} where we verify
special cases of Question~\ref{tach03}\eqref{tach03item1}.
The primary result of that section
is the following theorem whose proof is in~\ref{ABS2.1}. 
Recall that $R$ is \emph{generically Gorenstein} if,
for each $\p\in\ass(R)$, the ring $R_{\p}$ is Gorenstein.

\begin{thm} \label{intthmC}
Let $R$ be a local ring that is generically Gorenstein
and admits a dualizing complex $D$ such that $\inf(D)=0$.
Fix  semidualizing $R$-complexes $B$ and $C$ such that $\inf(B)=0=\inf(C)$.
Assume that $B$ is Cohen-Macaulay and $\sup(C)=0$.
\begin{enumerate}[\quad\rm(a)]
\item \label{ABS2.1item1}
If  $\ext^n_R(B,\rhom_R(B,D))=0$
for $n=1,\ldots,\dim(R)$, then $B\simeq R$.
\item \label{ABS2.1item0}
If  $\ext^n_R(\rhom_R(C,D),C)=0$
for $n=1,\ldots,\dim(R)$, then $C\simeq D$.
\end{enumerate}
\end{thm}

Note that the special
case $B=D$  in part~\eqref{ABS2.1item1} 
(or $C=R$ in part~\eqref{ABS2.1item0})
is exactly~\cite[(2.1)]{avramov:edcrcvct}.
Of course, we do not extend all of the special cases covered in~\cite{avramov:edcrcvct}
to the semidualizing arena.
Rather, we prove a few  results of this type in order  to illustrate the natural parallels
between the two contexts.

Unlike much of the existing literature on the subject,
most of this paper is devoted to the study of semidualizing complexes
over nonlocal rings.
In a sense, this makes it a natural companion
to~\cite{frankild:rrhffd}. 
However, it should be noted that many of our results
are new even in the local case.

\section{Complexes}\label{secintro}

Throughout this paper $R$ is a commutative noetherian ring.

\begin{defn} \label{notn01}
We index $R$-complexes homologically
$$X =\cdots\xra{\partial^X_{n+1}}X_n\xra{\partial^X_n}
X_{n-1}\xra{\partial^X_{n-1}}\cdots$$
and the \emph{infimum}, \emph{supremum}, and \emph{amplitude} of 
an $R$-complex $X$ are
\begin{gather*}
\inf(X)=\inf\{i\in\zz\mid\HH_i(X)\neq 0\} \qquad
\sup(X)=\sup\{i\in\zz\mid\HH_i(X)\neq 0\} \\
\amp(X)=\sup(X)-\inf(X).
\end{gather*}
The complex $X$ is \emph{homologically bounded} if $\amp(X)<\infty$;
it is \emph{degreewise homologically finite} if each
$R$-module $\HH_n(X)$ is finitely generated;
and it is \emph{homologically finite} if the
$R$-module $\oplus_{n\in\zz}\HH_n(X)$ is finitely generated.

For each integer $i$,
the $i$th \emph{suspension} (or \emph{shift}) of
a complex $X$, denoted $\shift^i X$, is the complex with
$(\shift^i X)_n:=X_{n-i}$ and $\partial_n^{\shift^i X}:=(-1)^i\partial_{n-i}^X$.
The notation $\shift X$ is short for $\shift^1 X$.
The \emph{projective dimension},
\emph{flat dimension} and  \emph{injective dimension}
of $X$ are denoted $\pd_R(X)$, $\fd_R(X)$ and $\id_R(X)$, respectively;
see~\cite{avramov:hdouc}.
\end{defn}

\begin{defn} \label{notn02}
We work in the derived category $\D(R)$.
References 
on the subject include~\cite{gelfand:moha,hartshorne:rad,verdier:cd,verdier:1}.
The category 
$\catdb(R)$ is the full subcategory of $\D(R)$
consisting of homologically bounded $R$-complexes.
Given two $R$-complexes $X$ and $Y$,
the derived homomorphism and tensor product complexes
are denoted $\rhom_R(X,Y)$ and $X\lotimes_R Y$, respectively.
For each integer $n$, set
$$\ext^n_R(X,Y):=\HH_{-n}(\rhom_R(X,Y)
\qquad
\text{and}
\qquad
\tor_n^R(X,Y):=\HH_n(X\lotimes_R Y).$$
Isomorphisms in $\D(R)$ are identified by the symbol $\simeq$,
and isomorphisms up to shift are identified by $\sim$.

The \emph{support} and \emph{dimension} of $X$ are, respectively,
\begin{align*}
\supp_R(X)
&=\{\p\in\spec(R)\mid X_{\p}\not\simeq 0\}
=\cup_n\supp_R(\HH_n(X)) \\
\dim_R(X)
&=\sup\{\dim(R/\p)-\inf(X_{\p})\mid\p\in\supp_R(X)\}. 
\end{align*}
\end{defn}

\begin{defn} \label{notn02'}
Assume that $(R,\m,k)$ is a local ring, and let $X$ be an $R$-complex.
The \emph{depth} and 
\emph{Cohen-Macaulay defect} of $X$ are, respectively
\begin{align*}
\depth_R(X)
&=-\sup(\rhom_R(k,X)) \\
\cmd_R(X)
&=\dim_R(X)-\depth_R(X).
\end{align*}
When $X$ is homologically finite, we have
$\cmd_R(X)\geq 0$ by~\cite[(2.8), (3.9)]{foxby:bcfm},
and $X$ is \emph{Cohen-Macaulay} if 
$\cmd_R(X)=0$.
\end{defn}

\begin{fact} \label{ext03}
For $R$-complexes $X$ and $Y$, 
the following are from~\cite[(2.1)]{foxby:ibcahtm}
\begin{align*}
\sup(\rhom_R(X,Y))&\leq\sup(Y)-\inf(X)\\
\inf(X\lotimes_RY)&\geq\inf(X)+\inf(Y)\\
\ext^{\inf(X)-\sup(Y)}_R(X,Y)&\cong\Hom_R(\HH_{\inf(X)}(X),\HH_{\sup(Y)}(Y)) \\
\tor^R_{\inf(X)+\inf(Y)}(X,Y)&\cong\HH_{\inf(X)}(X)\otimes_R\HH_{\inf(Y)}(Y).
\end{align*}
Assume that $R$ is local and that $X$ and $Y$ are 
degreewise homologically  finite such that 
$\inf(X),\inf(Y)>-\infty$.
Nakayama's Lemma and
the previous display imply
$\inf(X\lotimes_R Y)=\inf(X)+\inf(Y)$.
\end{fact}

\begin{fact} \label{ext03'}
Assume that $R$ is local and that $X$, $Y$ and $Z$ are 
degreewise homologically  finite $R$-complexes such that
$\inf(X),\inf(Y)>-\infty$ and $\sup(Z)<\infty$.
Using~\cite[(1.5.3)]{avramov:rhafgd} 
we see that $\pd_R(X\lotimes_R Y)<\infty$ if and only if
$\pd_R(X)<\infty$ and $\pd_R(Y)<\infty$, and
$\id_R(\rhom_R(X, Z))<\infty$ if and only if
$\pd_R(X)<\infty$ and $\id_R(Z)<\infty$.
\end{fact}

We shall have several occasions to use the following isomorphisms from~\cite[(4.4)]{avramov:hdouc}.

\begin{notationdefinition} \label{basics03}
Let  $X$, $Y$ and $Z$ be $R$-complexes.  Assume that $X$
is degreewise homologically finite and $\inf(X)>-\infty$.

The natural
\emph{tensor-evaluation morphism}
$$\omega_{XYZ}\colon\rhom_R(X,Y)\lotimes_RZ\to\rhom_R(X,Y\lotimes_RZ)$$
is an isomorphism when $\sup(Y)<\infty$ and either 
$\pd_R(X)<\infty$ or $\pd_R(Z)<\infty$.

The natural
\emph{Hom-evaluation morphism}
$$\theta_{XYZ}\colon X\lotimes_R\rhom_R(Y,Z)\to\rhom_{R}(\rhom_R(X,Y),Z)$$
is an isomorphism when $Y\in\D_{\mathrm{b}}(R)$ and either 
$\pd_R(X)<\infty$ or $\id_R(Z)<\infty$.
\end{notationdefinition}

Semidualizing complexes, defined next, are our main objects of study.

\begin{defn} \label{notn06}
A homologically finite $R$-complex $C$
is \emph{semidualizing} if the natural homothety morphism
$\chi^R_C\colon R\to \rhom_R(C,C)$
is an isomorphism in $\D(R)$.
An $R$-complex $D$ is \emph{dualizing} if it is semidualizing
and $\id_R(D)<\infty$.
Let $\s(R)$ denote the set of isomorphism classes of semidualizing
$R$-complexes.
\end{defn}

\begin{disc} \label{disc0601}
In some of the literature, 
the set of \emph{shift}-isomorphism classes of semidualizing $R$-complexes
is denoted $\mathfrak{S}(R)$.  
The notation $\s(R)$  is meant to evoke the notation $\mathfrak{S}(R)$ 
while at the same time distinguishing between the two notations.
\end{disc}

We include the following properties for ease of reference.

\begin{props} 
Let $C$ be a semidualizing $R$-complex.

\begin{subprops}
\label{sdc01item0}
The $R$-module $R$
is $R$-semidualizing. 
When $R$ is local, 
we have $\pd_R(C)<\infty$
if and only if $C\sim R$ by~\cite[(8.1)]{christensen:scatac}.  
\end{subprops}

\begin{subprops}
\label{dual03item2}
If $R$ is Gorenstein and local, then $C\sim R$.  
Conversely, if $R$ is dualizing for $R$, then $R$ is Gorenstein.
See~\cite[(8.6)]{christensen:scatac} and~\cite[(V.9)]{hartshorne:rad}.
\end{subprops}

\begin{subprops} \label{sp02}
If $X$ is a homologically finite $R$-complex, then 
$X$ is semidualizing for $R$ if and only if
$X_{\m}$ is semidualizing for $R_{\m}$ for each maximal
(equivalently, for each prime) ideal $\m\subset R$;
see~\cite[(2.3)]{frankild:rrhffd}.
When $R$ is local and $C$ is semidualizing with
$s=\sup(C)$, if $\p\in\ass_R(\HH_s(C))$,
then $\inf(C_{\p})=s$ by~\cite[(A.7)]{christensen:scatac}.
\end{subprops}

\begin{subprops} \label{sp03}
Let $\vf\colon R\to S$ be a local homomorphism
of finite flat dimension, and fix semidualizing
$R$-complexes $B,C$. The complex $S\lotimes_R C$ is
semidualizing for $S$, and 
$\amp(S\lotimes_R C)=\amp(C)$
by~\cite[(5.7)]{christensen:scatac}.
If $S\lotimes_R C$ is
dualizing for $S$, then
$C$ is dualizing for $R$ by~\cite[(4.2), (5.1)]{avramov:lgh}.
Conversely, if $C$ is dualizing for $R$ and 
$\vf$ is surjective with kernel generated by an $R$-sequence,
then $S\lotimes_R C$ is
dualizing for $S$, 
and $\comp{R}\lotimes_R C$ is
dualizing for $\comp{R}$ 
by~\cite[(4.2), (4.3), (5.1)]{avramov:lgh}.
Finally, if $S\lotimes_R B\simeq S\lotimes_R C$ in $\D(S)$,
then~\cite[(1.10)]{frankild:rrhffd} implies that
$B\simeq C$ in $\D(R)$.
\end{subprops}

\begin{subprops} \label{sp01}
Let $\alpha\colon X\to Y$ be a morphism between degreewise
homologically finite $R$-complexes.
Assuming $\inf(X),\inf(Y)>-\infty$, if $\rhom(\alpha,C)$ is an
isomorphism in $\D(R)$, then so is $\alpha$.
Dually, assuming $\sup(X),\sup(Y)<\infty$, if $\rhom(C,\alpha)$ is an
isomorphism in $\D(R)$, then so is $\alpha$.
These follow from~\cite[(1.2.3.b)]{avramov:rhafgd}
and~\cite[(A.8.11), (A.8.13)]{christensen:gd}
as the isomorphism $R\simeq\rhom_R(C,C)$
implies $\supp_R(C)=\spec(R)$.
\end{subprops}

\begin{subprops} \label{dual04''}
If $C$ is a module, then an element of $R$ is
$C$-regular if and only if it is $R$-regular
as the isomorphism $R\cong\hom_R(C,C)$
implies $\ass_R(C)=\ass(R)$.
\end{subprops}

\begin{subprops} \label{dual04'}
When $R$ is local, there are inequalities
\begin{align*}
\max\{\amp(C),\cmd_R(C)\}&\leq\cmd(R)\leq\amp(C)+\cmd_R(C)
\end{align*}
\label{sdc01item2}
by~\cite[(3.4)]{christensen:scatac}.
In particular, if $R$ is Cohen-Macaulay and local, then $\amp(C)=0$.
\end{subprops}
\end{props}

The  next definition is from~\cite{christensen:scatac}
and~\cite{hartshorne:rad} and 
will be used primarily to compare semidualizing complexes.

\begin{defn} \label{notn08}
Let $C$ be a semidualizing $R$-complex.
A homologically finite $R$-complex $X$ is \emph{$C$-reflexive} 
when it satisfies the 
following:
\begin{enumerate}
\item \label{notn08item2}
$\ext^n_R(X,C)=0$ for $n\gg 0$, and
\item \label{notn08item3}
the natural biduality morphism $\delta^C_X\colon X\to\rhom_R(\rhom_R(X,C),C)$
is an isomorphism in $\D(R)$.
\end{enumerate}
\end{defn}

The following properties are frequently used in the sequel.

\begin{props} 
Let $B$ and $C$ be semidualizing $R$-complexes.

\begin{subprops}
\label{sdc02item1}
Each homologically finite $R$-complex of finite projective
dimension is $C$-reflexive by~\cite[(3.11)]{frankild:rrhffd}.
\end{subprops}

\begin{subprops}
\label{sdc02item3}
Let $D$ be a dualizing $R$-complex, and $X$ a homologically finite 
$R$-complex. The complex $X$ is $D$-reflexive by~\cite[(V.2.1)]{hartshorne:rad},
and $X$ is semidualizing
if and only if $\rhom_R(X,D)$ is semidualizing by~\cite[(2.12)]{christensen:scatac}
and~\eqref{sp02}.
The natural evaluation morphism
$C\lotimes_R \rhom_R(C,D)\to D$
is an isomorphism in $\D(R)$ by~\cite[(3.1.b)]{frankild:rrhffd}.

Assume that $R$ is local and $\inf(C)=0=\inf(D)$.
With the isomorphism from the previous paragraph,
Fact~\ref{ext03} yields the second equality in the next sequence 
$$0=\inf(D)=\inf(C)+\inf(\rhom_R(C,D))=\inf(\rhom_R(C,D))
$$
and~\cite[(3.14)]{foxby:bcfm} implies
$\amp(\rhom_R(X,D))=\cmd(X)$.
\end{subprops}

\begin{subprops}
\label{sdc02item5}
If $X$ is a homologically finite $R$-complex,
then $X$ is $C$-reflexive if and only if 
$\ext_R^i(X,C)=0$ for $i\gg 0$ and
$X_{\m}$ is $C_{\m}$-reflexive
for each maximal (equivalently, for each prime)
ideal $\m\subset R$; see~\cite[(2.8)]{frankild:rrhffd}. 
\end{subprops}

\begin{subprops}
\label{sdc02item2}
If  $R$ is local, then 
$B$ is $C$-reflexive and $C$ is $B$-reflexive
if and only if $B\sim C$ by~\cite[(5.5)]{takahashi:hiatsb}.
\end{subprops}

\begin{subprops}
\label{sdc02item4}
When $B\lotimes_R C$ is also semidualizing
we know from~\cite[(3.1.c)]{frankild:rrhffd} that
$C$ is
$B\lotimes_RC$-reflexive
and furthermore 
$\rhom_R(C,B\lotimes_RC)\simeq B$.
\end{subprops}
\end{props}

The following categories, known collectively as
``Foxby classes'', 
were introduced  by Foxby~\cite{foxby:gdcmr},
Avramov and Foxby~\cite{avramov:rhafgd},
and Christensen~\cite{christensen:scatac}.

\begin{defn} \label{bass01}
Let $C$ be a semidualizing $R$-complex.
The \emph{Auslander class} associated to $C$, 
denoted $\catac(R)$, is the full subcategory of $\D(R)$
consisting of all complexes $X$ satisfying the following conditions:
\begin{enumerate}[\quad(1)]
\item \label{bass01item1}
$X$ and $C\lotimes_R X$ are homologically bounded, and
\item \label{bass01item2}
the natural morphism $\gamma^C_X\colon X\to\rhom_R(C,C\lotimes_R X)$
is an isomorphism.
\end{enumerate}
The \emph{Bass class}  associated to $C$, 
denoted $\catbc(R)$, is the full subcategory of $\D(R)$
consisting of all complexes $X$ satisfying the following conditions:
\begin{enumerate}[\quad(1)]
\item \label{bass01item3}
$X$ and $\rhom_R(C,X)$ are homologically bounded, and
\item \label{bass01item4}
the evaluation morphism $\xi^C_X\colon C\lotimes_R \rhom_R(C,X)\to X$
is an isomorphism.
\end{enumerate}
\end{defn}

\begin{fact} \label{bass02}
Let  $C$ be a semidualizing $R$-complex.
From~\cite[(4.6)]{christensen:scatac}
we conclude that a homologically bounded complex $X$ is in 
$\catac(R)$ if and only if $C\lotimes_R X$ is in $\catbc(R)$,
and dually, a homologically bounded complex $X$ is in $\catbc(R)$
if and only if $\rhom_R(C,X)$ is in $\catac(R)$.
Each $R$-complex of finite injective dimension
is in $\catbc(R)$, and each $R$-complex of finite projective dimension
is in $\catac(R)$ by~\cite[(4.4)]{christensen:scatac}.
\end{fact}

The complete intersection dimensions used in this paper
were defined for modules by
Avramov, Gasharov and Peeva~\cite{avramov:cid} and 
Sahandi, Sharif and Yassemi~\cite{sahandi:hfd},
and then for complexes by Sather-Wagstaff~\cite{sather:cidfc,sather:cidc}.  
We start with the definitions, first over a local ring and then in general.

\begin{defn} \label{cidim01}
Assume that $R$ is local.
A \emph{quasi-deformation} of  $R$ is a diagram of local
ring homomorphisms $R\xra{\vf} R'\xla{\rho} Q$ in which
$\vf$ is flat and $\rho$ is surjective with $\ker(\rho)$ generated
by a $Q$-regular sequence.

For each homologically finite $R$-complex $X$,
define
the 
\emph{complete intersection dimension} and
\emph{complete intersection injective dimension} of $X$ as
follows.
\begin{align*}
\cidim_R(X)
&:=\inf\left\{\pd_Q(R'\lotimes_R X)-\pd_Q(R')\left| \text{
\begin{tabular}{@{}c@{}}
$R\to R'\from Q$ is a \\ quasi-deformation 
\end{tabular}
}\!\!\!\right. \right\} \\
\ciid_R(X)
&:=\inf\left\{\id_Q(R'\lotimes_R X)-\pd_Q(R')
\left| \text{
\begin{tabular}{@{}c@{}}
$R\to R'\from Q$ is a \\ quasi-deformation 
\end{tabular}
}\!\!\!\right. \right\}
\end{align*}
\end{defn}

\begin{defn} \label{cidim01'}
For each homologically finite $R$-complex $X$,
define
the 
\emph{complete intersection dimension} and
\emph{complete intersection injective dimension} of $X$ as
\begin{align*}
\cidim_R(X)
&:=\sup\left\{\cidim_{R_{\m}}(X_{\m}) \mid \text{$\m\subset R$ is a maximal ideal}
\right\} \\
\ciid_R(X)
&:=\sup\left\{\ciid_{R_{\m}}(X_{\m}) \mid \text{$\m\subset R$ is a maximal ideal}
\right\}.
\end{align*}
\end{defn}

\begin{fact} \label{cidim02}
Let $X$ be a homologically finite $R$-complex. 
If $\pd_R(X)<\infty$, then $\cidim_R(X)<\infty$ by~\cite[(3.3)]{sather:cidc}.
If $R$ is local and 
$\id_R(X)<\infty$, then the trivial quasi-deformation 
$R\to R\from R$ yields
$\ciid_R(X)<\infty$;
see Lemma~\ref{solve01} for the nonlocal situation.

Assume that $R$ is local.
One checks readily that $\cidim_R(X)$ is finite
if and only if $\cidim_{\comp R}(\comp{R}\lotimes_RX)$ is finite. 
If $X$ has finite length homology, then $\ciid_R(X)$ is finite
if and only if $\ciid_{\comp R}(\comp{R}\lotimes_RX)$ is finite
by~\cite[(3.7)]{sather:cidfc}.
\end{fact}

We document several useful but independent
facts about quasi-deformations and complete intersection dimensions
in Appendix~\ref{sec07}.

\section{Detecting Reflexivity} \label{sec03}

The results of this section are of the first type discussed in the introduction:
Assuming $\ext^n_R(B,C)=0$ for $n\gg 0$ \emph{and a bit more},
we show that $B$ is $C$-reflexive. 
We begin, though, by discussing the link 
between~\eqref{ABSq} and Question~\ref{tach03}.

\begin{disc} \label{sdc03}
Assume that the answer to Question~\ref{tach03}\eqref{tach03item1} is ``yes''.
Then the answer to the question in~\eqref{ABSq} is also ``yes''.
Indeed, let $D$
be a dualizing $R$-complex such that $\inf(D)=0$,
and assume that  $\ext^n_R(D,R)=0$ for $(\dim(R)+1)$
consecutive values of  $n\geq 1$.
The $R$-complexes $R$ and $D$ are semidualizing,
and $\sup(R)=0=\inf(D)$.
Thus, the affirmative answer to Question~\ref{tach03}\eqref{tach03item1} implies that
$D$ is $R$-reflexive.  From~\eqref{sdc02item3}
we know that $R$ is $D$-reflexive, and so~\eqref{sdc02item2} implies
$R\sim D$.
Using~\eqref{dual03item2} we conclude that $R$ is
Gorenstein, as desired.

Similarly, if the answer to Question~\ref{tach03}\eqref{tach03item2} is ``yes'',
then this would establish the conjecture in~\eqref{ABSq}.
\end{disc}

Before proceeding, we briefly discuss the necessity of the
semidualizing hypothesis in
Question~\ref{tach03}.

\begin{disc} \label{sega01}
If
the
semidualizing-hypothesis is removed from Question~\ref{tach03},
then the answer to each of the resulting questions is ``no''.
Specifically, Jorgensen and \c{S}ega~\cite[Theorem]{jorgensen:itrcm}
exhibit a local ring $R$ and a family $\{M_s\}_{2\leq s\leq \infty}$
of $R$-modules such that $\ext^n_R(M_s,R)=0$ for $i=1,\ldots,s$
and $M_s$ is not $R$-reflexive, for each $s$.
\end{disc}

When $R$ is local,
the forward implication in the following result is
in~\cite[(2.11)]{christensen:scatac},
but the proof of~\cite[(2.11)]{christensen:scatac}
makes no use of the local hypothesis. 
For the reverse implication, argue as in~\cite[(2.1.10)]{christensen:gd}.

\begin{lem}\label{reflexive1}
Let $X$ and $C$ be
homologically finite $R$-complexes  with $C$ semidualizing. 
The complex $X$ is $C$-reflexive if and only if $\rhom_R(X,C)$ is 
$C$-reflexive.
\qed
\end{lem}

The following diagram will be used in the next two proofs.

\begin{disc} \label{diagram01}
Let $X$ and $C$ be 
homologically finite $R$-complexes  with $C$ semidualizing. 
There is a commutative diagram of morphisms of complexes
$$\xymatrix{
R \ar[rr]^-{\chi^R_{\rhom_R(X,C)}} \ar[d]_{\chi^R_X} 
&& \rhom_R(\rhom_R(X,C),\rhom_R(X,C)) \ar[d]^{\simeq} \\
\rhom_R(X,X) \ar[rr]^-{\rhom_R(X,\delta^C_X)}
&& \rhom_R(X,\rhom_R(\rhom_R(X,C),C))
}
$$
wherein the unspecified arrow is the ``swap''
isomorphism from~\cite[(A.2.9)]{christensen:gd}.
\end{disc}

Before proving Theorem~\ref{intthmA}, we present a similar
result, where only one of the complexes is assumed
to be semidualizing. When $R$ is local,
the implication
\eqref{semidual1item1}$\implies$\eqref{semidual1item2} 
is in~\cite[(2.11)]{christensen:scatac},
but the proof of ~\cite[(2.11)]{christensen:scatac}
makes no use of the local hypothesis. 
Also, the symmetry of the conditions in our result suggest a fourth
condition, namely, 
``$X$ is semidualizing and $\rhom_R(X,C)$ is $C$-reflexive'';
this is
shown to be equivalent in Theorem~\ref{intthmA}.

\begin{prop}\label{semidual1}
Let $X$ and $C$ be
homologically finite $R$-complexes  with $C$ semidualizing. 
The following conditions are equivalent:
\begin{enumerate}[\quad\rm(i)]
\item \label{semidual1item1}
$X$ is semidualizing and $C$-reflexive;
\item \label{semidual1item2}
$\rhom_R(X,C)$  is  semidualizing and $C$-reflexive;
\item \label{semidual1item4}
$\rhom_R(X,C)$ is semidualizing and $X$  is $C$-reflexive.
\end{enumerate}
\end{prop}

\begin{proof} 
As noted above, the implication 
\eqref{semidual1item1}$\implies$\eqref{semidual1item2}
is proved as  in~\cite[(2.11)]{christensen:scatac}. 
Also, the equivalence
\eqref{semidual1item2}$\iff$\eqref{semidual1item4}
is from Lemma~\ref{reflexive1}.

\eqref{semidual1item4}$\implies$\eqref{semidual1item1}.
Assume that $\rhom_R(X,C)$ is semidualizing and $X$  is $C$-reflexive.
Then the morphisms $\chi^R_{\rhom_R(X,C)}$ and $\delta^C_X$
are isomorphisms.
Hence, the morphism $\rhom_R(X,\delta^C_X)$ is an isomorphism, and so
the diagram in Remark~\ref{diagram01} shows that $\chi^R_X$ is also an isomorphism.
By definition, we conclude that $X$ is semidualizing.
\end{proof}

\begin{para} \label{semidual2}
\emph{Proof of Theorem~\ref{intthmA}.}
\eqref{semidual2item1}$\iff$\eqref{semidual2item3}.  This
is from Lemma~\ref{reflexive1}. 

\eqref{semidual2item1}$\iff$\eqref{semidual2item2}. 
The conditions~\eqref{semidual2item1} and~\eqref{semidual2item2}
each imply that $\rhom_R(B,C)$ is homologically bounded.
So, it remains to assume
that $\rhom_R(B,C)$ is homologically bounded and
show that $\chi_{\rhom_R(B,C)}^R$ and $\delta_B^C$
are isomorphisms simultaneously.

As $B$ is semidualizing,
the morphism $\chi_B^R$ is an isomorphism.  Hence
the diagram from Remark~\ref{diagram01} with $X=B$
shows that  $\chi_{\rhom_R(B,C)}^R$ 
and $\rhom_R(B,\delta_B^C)$  are isomorphisms simultaneously. 
From~\eqref{sp01}, we know that
$\rhom_R(B,\delta_B^C)$ is an isomorphism
if and only if $\delta_B^C$ is so, and hence the desired
equivalence.

\eqref{semidual2item2}$\iff$\eqref{semidual2item4}.
Consider the following commutative diagram of morphisms of 
complexes wherein the unspecified isomorphism is a combination 
of Hom-tensor adjointness and commutativity of tensor product.
$$\xymatrix{
R \ar[rr]^-{\chi^R_{\rhom_R(B,C)}} \ar[d]_{\chi^R_C}^{\simeq} 
&& \rhom_R(\rhom_R(B,C),\rhom_R(B,C)) \ar[d]^{\simeq} \\
\rhom_R(C,C) \ar[rr]^-{\rhom_R(\xi^B_C,C)}
&& \rhom_R(B\lotimes_R\rhom_R(B,C),C)
}
$$
Using this diagram with~\eqref{sp01},
the desired equivalence is established as in the proof of
\eqref{semidual2item1}$\iff$\eqref{semidual2item2}.  \qed
\end{para}

Next we prove versions of Proposition~\ref{semidual1}
and Theorem~\ref{intthmA} for Auslander classes. 

\begin{prop}\label{semidual5}
Let $X$ and $C$  be
homologically finite $R$-complexes with $C$ semidualizing. 
The following conditions are equivalent:
\begin{enumerate}[\quad\rm(i)]
\item \label{semidual5item1}
$X$ is semidualizing and in $\catac(R)$;
\item \label{semidual5item2}
$X$ and $C\lotimes_RX$ are semidualizing;
\item \label{semidual5item3}
$X$  is in $\catac(R)$ and $C\lotimes_RX$ is semidualizing;
\item \label{semidual5item4}
$C\lotimes_RX$ is semidualizing and $C$ is $C\lotimes_RX$-reflexive;
\item \label{semidual5item5}
$C\lotimes_RX$ is semidualizing and $X$ is $C\lotimes_RX$-reflexive.
\end{enumerate}
\end{prop}

\begin{proof} 
Consider the following commutative diagram of morphisms
of complexes 
$$\xymatrix{
R \ar[rr]^-{\chi^R_{C\lotimes_R X}} \ar[d]_{\chi^R_X} 
&& \rhom_R(C\lotimes_R X,C\lotimes_R X) \ar[d]^{\simeq} \\
\rhom_R(X,X) \ar[rr]^-{\rhom_R(X,\gamma^C_X)}
&& \rhom_R(X,\rhom_R(C,C\lotimes_R X))
}
$$
wherein the unspecified isomorphism is a combination 
of Hom-tensor adjointness and commutativity of tensor product.

\eqref{semidual5item1}$\iff$\eqref{semidual5item2}$\iff$\eqref{semidual5item3}.
Use the above diagram
as in~\eqref{semidual2}. 

\eqref{semidual5item3}$\implies$\eqref{semidual5item4}
and \eqref{semidual5item3}$\implies$\eqref{semidual5item5}.
Assume that $X$  is in $\catac(R)$ and $C\lotimes_RX$ is semidualizing.
Using the above diagram, we  see that $X$ is semidualizing,
and~\eqref{sdc02item4} implies that 
$C$ and $X$ are
$C\lotimes_RX$-reflexive.

\eqref{semidual5item4}$\implies$\eqref{semidual5item3}.
Assume that $C\lotimes_RX$ is semidualizing and $C$ is $C\lotimes_RX$-reflexive.
Theorem~\ref{intthmA} implies $C\lotimes_R X\in\catbc(R)$, and so
Fact~\ref{bass02} implies $X\in\catac(R)$.  

\eqref{semidual5item5}$\implies$\eqref{semidual5item4}.
Assume that
$C\lotimes_RX$ is semidualizing and $X$ is $C\lotimes_RX$-reflexive.
The morphism 
$\gamma^X_C\colon C\to\rhom_R(X,X\lotimes_R C)$
is locally an isomorphism 
by a result of Gerko~\cite[(3.5)]{gerko:sdc}, and hence 
$\gamma^X_C$ is an isomorphism.
Since $X$ is $C\lotimes_R X$-reflexive, we conclude from
Lemma~\ref{reflexive1} that the following complex is also $C\lotimes_R X$-reflexive
\begin{xxalignat}{3}
  &{\hphantom{\square}}& \rhom_R(X,C\lotimes_R X)
  &\simeq \rhom_R(X,X\lotimes_R C)\simeq C. &&\qedhere
\end{xxalignat}
\end{proof}

\begin{cor}\label{semidual6} 
Let $B$ and $C$ be
semidualizing $R$-complexes.
The following conditions are equivalent:
\begin{enumerate}[\quad\rm(i)]
\item \label{semidual6item1}
$C\lotimes_R B$  is  semidualizing;
\item \label{semidual6item2}
$B$ is in $\catac(R)$;
\item \label{semidual6item3}
$C$ is in $\catab(R)$.
\qed
\end{enumerate}
\end{cor}

Observe that the implication~\eqref{semidual6item2}$\implies$\eqref{semidual6item1}
in Corollary~\ref{semidual6} has the following form:
Assuming $\tor_n^R(B,C)=0$ for $n\gg 0$ \emph{and a bit more},
we conclude that $B\lotimes_RC$ is semidualizing.
In light of Question~\ref{tach03} and Theorem~\ref{intthmA},
this motivates the next question.
Regarding the hypotheses of Question~\ref{tpq02}\eqref{tpq02item1},
Proposition~\ref{tor01} 
gives a partial justification of
the range of Tor-vanishing, 
and 
Example~\ref{local01} shows why $R$ must be local.

\begin{question} \label{tpq02}
Let  $B$ and $C$ be
semidualizing $R$-complexes.
\begin{enumerate}[\quad(a)]
\item \label{tpq02item1}
If $\tor^R_n(B,C)=0$ 
for $(2\dim(R)+1)$ consecutive values of $n>\inf(B)+\inf(C)$
and $R$ is local,
must $B\lotimes_RC$ be semidualizing?
\item \label{tpq02item2}
If $\tor^R_n(B,C)=0$
for all  $n>\inf(B)+\inf(C)$,
must $B\lotimes_RC$ be semidualizing?
\item \label{tpq02item3}
If $\tor^R_n(B,C)=0$ for $n\gg 0$,
must $B\lotimes_RC$ be semidualizing?
\end{enumerate}
\end{question}

We raise the next questions in light of Propositions~\ref{semidual1}
and~\ref{semidual5}.

\begin{question} \label{tpq01}
Let $X$ and $C$ be
homologically finite $R$-complexes, and assume that $C$ is semidualizing. 
\begin{enumerate}[\quad\rm(a)]
\item \label{tpq01item1}
If $C\lotimes_RX$ is semidualizing, must $X$ also be
semidualizing?
\item \label{tpq01item2}
If $\rhom_R(C,X)$ is semidualizing, must $X$ also be
semidualizing?
\item \label{tpq01item3}
If $\rhom_R(X,C)$ is semidualizing, must $X$ also be
semidualizing?
\end{enumerate}
\end{question}

\begin{disc} \label{disc0301}
As with Question~\ref{tach03}, if we assume
more in Question~\ref{tpq01}, then we have
an affirmative answer.  For instance, if $\cidim_R(X)$ is finite
and $C\lotimes_RX$ is semidualizing, then $X$ is also 
semidualizing. Indeed, using~\eqref{sp02} we may assume
without loss of generality that $R$ is local.  The finiteness
of $\cidim_R(X)$ implies $X\in\catac(R)$ by~\cite[(5.1.a)]{sather:cidfc}.
Hence, Proposition~\ref{semidual5} implies 
that $X$ is semidualizing. 

Similarly, if  either (1) $\cidim_R(X)<\infty$
and $\rhom_R(X,C)$ is semidualizing or (2) $\ciid_R(X)<\infty$
and $\rhom_R(C,X)$ is semidualizing, then $X$ is 
semidualizing.
\end{disc}

The final result of
this section shows that Questions~\ref{tach03}\eqref{tach03item3}
and~\ref{tpq02}\eqref{tpq02item3} are equivalent when $R$ admits a
dualizing complex.

\begin{prop} \label{tpq03}
Assume that $R$  admits a dualizing complex $D$.
The following conditions are equivalent:
\begin{enumerate}[\quad\rm(i)]
\item \label{tpq03item1}
For all semidualizing
$R$-complexes $B$ and $C$, if $\ext_R^n(B,C)=0$ for $n\gg 0$,
then $\rhom_R(B,C)$ is semidualizing;
\item \label{tpq03item2}
For all semidualizing
$R$-complexes $B$ and $C$, if $\tor^R_n(B,C)=0$ for $n\gg 0$,
then $B\lotimes_R C$ is semidualizing.
\end{enumerate}
\end{prop}

\begin{proof}
\eqref{tpq03item1}$\implies$\eqref{tpq03item2}.
Assume $\tor^R_n(B,C)=0$ for $n\gg 0$.  This means that the complex
$B\lotimes_R C$ is homologically finite.  It follows that the same
is true for the complexes in the next display
where the isomorphism is Hom-tensor adjointness
$$\rhom_R(B\lotimes_R C,D)\simeq\rhom_R(B,\rhom_R(C,D)).$$
The homological finiteness of the second complex says
$\ext_R^n(B,\rhom_R(C,D))=0$ for $n\gg 0$.
Because the complexes
$B$ and $\rhom_R(C,D)$ are semidualizing, condition~\eqref{tpq03item1}
implies that $\rhom_R(B,\rhom_R(C,D))$ is semidualizing.
The displayed isomorphism implies that $\rhom_R(B\lotimes_R C,D)$ is semidualizing,
and so $B\lotimes_R C$ is semidualizing by~\eqref{sdc02item3}.

\eqref{tpq03item2}$\implies$\eqref{tpq03item1}.
The proof is similar to the previous paragraph using the 
following Hom-evaluation isomorphism~\eqref{basics03}
in place of adjointness
\begin{xxalignat}{3}
  &{\hphantom{\square}}& \rhom_R(\rhom_R(B,C),D)
  &\simeq B\lotimes_R\rhom_R(C,D). &&\qedhere
\end{xxalignat}
\end{proof}

\section{Derived Picard Group Action and Associated Relation} \label{sec05}

The goal of this section is to discuss a 
relation $\approx$ that is better suited than $\sim$ for studying
semidualizing complexes over nonlocal rings.  
As motivation for our first definition, recall that a \emph{line bundle} on $R$
is a finitely generated 
locally free (i.e.,
projective) $R$-module of rank 1.  
The \emph{Picard group} of $R$, denoted $\Pic(R)$,
is the set of
isomorphism classes of  line bundles on $R$.  Given a 
line bundle $P$, its class in $\Pic(R)$ is 
denoted $[P]$.
As the name ``Picard group''
suggests, $\Pic(R)$ carries the structure of 
an abelian group with additive identity $[R]$ and operations
\[ [P]+[Q]=[P\otimes_R Q] \qquad\text{and}\qquad
[P]-[Q]=[\Hom_R(Q,P)]. \]
When $R$ is local, its Picard group is trivial, that is, we have $\Pic(R)=\{[R]\}$.

\begin{defn} \label{picard01}
A \emph{tilting $R$-complex} is
a homologically finite
$R$-complex $P$ of finite projective dimension
such that $P_{\p}\sim R_{\p}$ for each $\p\in\spec(R)$.
The \emph{derived Picard group} of $R$, denoted $\dpic(R)$,
is the set of
isomorphism classes of tilting
$R$-complexes,
and the class of a tilting $R$-complex $P$ in $\dpic(R)$ is 
denoted $[P]$.
For each $[P],[Q]\in\dpic(R)$, set
\[ [P]+[Q]=[P\lotimes_R Q] \qquad\text{and}\qquad
[P]-[Q]=[\rhom_R(Q,P)].\]
\end{defn}

\begin{disc} \label{disc0602}
Definition~\ref{picard01} is inspired by
the derived Picard group of Yekutieli~\cite{yekutieli:dcmedpgr}.  
Note that
our definition differs from Yekutieli's in that we do not assume that $R$ contains a field
and the definition does not depend on $R$ being an algebra.
\end{disc}

The following properties are verified using routine arguments.

\begin{properties} \label{picard02}
Let  $C$ be a semidualizing $R$-complex.
\begin{subprops} \label{picard02a}
The operations defined in~\eqref{picard01} endow $\dpic(R)$ with the structure of 
an abelian group with identity $[R]$. Tensor-evaluation~\eqref{basics03}
shows that $[\rhom_R(P,R)]$ is an additive inverse for
$[P]\in\dpic(R)$.
\end{subprops}
\begin{subprops} \label{picard02b}
Each line bundle on $R$ is a tilting complex, and this yields a 
well-defined injective abelian group homomorphism $\Pic(R)\hookrightarrow\dpic(R)$
given by $[P]\mapsto[P]$.
\end{subprops}
\begin{subprops} \label{picard02c}
Each tilting $R$-complex $P$ is homologically finite and locally semidualizing,
so~\eqref{sp02} implies $P$ is semidualizing. This yields a 
well-defined injective map $\dpic(R)\hookrightarrow\s(R)$
given by $[P]\mapsto[P]$.
\end{subprops}
\begin{subprops} \label{picard03}
If $R$ is local, then
$\dpic(R)=\{[\shift^nR]\mid n\in\zz\}\cong \zz$.
\end{subprops}
\end{properties}

The first result of this section shows that,
in order to verify that an  $R$-complex is tilting,
one need not check that it has finite projective dimension.

\begin{prop} \label{picard09}
A homologically finite $R$-complex $X$ is tilting if and only if
$X_{\m}\sim R_{\m}$ for each maximal ideal $\m\subset R$.
\end{prop}

\begin{proof}
For the nontrivial implication, assume 
$X_{\m}\sim R_{\m}$ for each maximal ideal $\m\subset R$.
It follows that
$X_{\p}\sim R_{\p}$ for each $\p\in\spec(R)$.
Because  $X$ is homologically finite,
this implies
$$\pd_{R_{\p}}(X_{\p})=\sup(X_{\p})\leq\sup(X).$$
From this we deduce the first inequality in the next sequence
$$\pd_R(X)
=\textstyle \sup_{\p}(\pd_{R_{\p}}(X_{\p}))\leq\sup(X)<\infty$$
while the equality is from~\cite[(5.3.P)]{avramov:hdouc}.
Thus $X$ is a tilting $R$-complex.
\end{proof}

Next, we prove a nonlocal version of~\cite[(8.1)]{christensen:scatac}.

\begin{lem} \label{tilt01}
Let $C$ be a semidualizing $R$-complex.
The following are equivalent:
\begin{enumerate}[\quad\rm(i)]
\item \label{tilt01a}
$C$ is a tilting $R$-complex;
\item \label{tilt01b}
$\catb_{C}(R)=\catdb(R)$;
\item \label{tilt01b'}
$R\in\catb_{C}(R)$;
\item \label{tilt01c}
$\cata_{C}(R)=\catdb(R)$;
\item \label{tilt01c'}
$E_R(R/\m)\in\cata_{C}(R)$ for each maximal ideal $\m\subset R$.
\end{enumerate}
\end{lem}

\begin{proof}
\eqref{tilt01a}$\implies$\eqref{tilt01b}.
Assume that $C$ is a tilting $R$-complex.
The containment $\catb_C(R)\subseteq\catdb(R)$
is by definition.
For the reverse containment,
fix a complex $X\in \catdb(R)$.  
Because $C$ and $X$ are both homologically bounded 
and $C$ has finite projective dimension we conclude that $\rhom_R(C,X)$
is homologically bounded. Consider the commutative diagram
$$\xymatrix{
C\lotimes_R\rhom_R(C,X) \ar[r]^-{\theta_{CCX}}_-{\simeq} \ar[d]_{\xi^C_X}
& \rhom_R(\rhom_R(C,C),X) \ar[d]_{\simeq}^{\rhom_R(\chi^R_C,X)} \\
X & \rhom_R(R,X) \ar[l]^-{\simeq}_-{\zeta_X}
}$$
wherein $\theta_{CCX}$ is the Hom-evaluation isomorphism~\eqref{basics03},
$\chi^R_C$ is the homothety isomorphism~\eqref{picard02c},
and $\zeta_X$ is the natural evaluation isomorphism.
The diagram shows that $\xi^C_X$ is an isomorphism, and so
$X\in\catb_C(R)$.

\eqref{tilt01b}$\implies$\eqref{tilt01b'}. This is immediate from the 
condition $R\in\catdb(R)$.

\eqref{tilt01b'}$\implies$\eqref{tilt01a}.
Assume that $R$ is in $\catbc(R)$.
Because $C$ is homologically finite, the condition $R\in\catbc(R)$
implies $R_{\m}\in\catb_{C_{\m}}(R_{\m})$ for each maximal ideal $\m\subset R$.
Hence, the local version of this 
result~\cite[(8.1)]{christensen:scatac} yields $C_{\m}\sim R_{\m}$
for each $\m$, and so $C$ is tilting by Proposition~\ref{picard09}.

The equivalences~\eqref{tilt01a}$\iff$\eqref{tilt01c}$\iff$\eqref{tilt01c'}
are established similarly.
\end{proof}

The following is
for use in the proofs of Theorem~\ref{intthmB}
and Proposition~\ref{sdc06}.

\begin{prop} \label{cidim03}
Let  $C$
be a semidualizing $R$-complex.
\begin{enumerate}[\quad\rm(a)]
\item  \label{cidim03item1}
$\cidim_R(C)<\infty$
if and only if 
$C$ is tilting.
\item  \label{cidim03item2}
$\ciid_R(C)<\infty$
if and only if $C$ is dualizing.
\end{enumerate}
\end{prop}

\begin{proof}
\eqref{cidim03item1}
One implication is in Fact~\ref{cidim02}.  
For the other implication, assume that
$\cidim_R(C)$ is finite.
It follows that
$\cidim_{R_{\m}}(C_{\m})<\infty$ for each maximal ideal $\m\subset R$.
If we can show that this implies $C_{\m}\sim R_{\m}$ for each $\m$,
then Proposition~\ref{picard09} will imply that $C$ is tilting.  Hence, we may replace
$R$ and $C$ by $R_{\m}$ and $C_{\m}$ in order to assume that $R$ is local.

Because $R$ is local and $\cidim_R(C)$ is finite, 
there exists a quasi-deformation
$R\to R'\from Q$ such that $\pd_Q(R'\lotimes_R C)$ is finite
and $Q$ is complete; see~\cite[Thm.\ F]{sather:cidfc}.  
Because $Q$ is complete, 
Lemma~\ref{lem0201} implies that there exists a semidualizing $Q$-complex
$N$ such that $R'\lotimes_R C\simeq R'\lotimes_Q N$.
The finiteness of $\pd_Q(R'\lotimes_Q N)=\pd_Q(R'\lotimes_R C)$
implies 
$\pd_Q(N)<\infty$ by Fact~\ref{ext03'}, and so
$N\sim Q$ by~\eqref{sdc01item0}. 
Hence, the isomorphism $C\sim R$ follows from
Lemma~\ref{lem0203}\eqref{lem0203a}.

\eqref{cidim03item2} 
For the nontrivial implication, assume $\ciid_R(C)<\infty$.
By Lemma~\ref{solve01}\eqref{solve01c}, it suffices to show that
$\id_{R_{\m}}(C_{\m})<\infty$ for each maximal ideal $\m\subset R$.
Hence, we may assume, as in the proof of part~\eqref{cidim03item1},
that $R$ is local with maximal ideal $\m$.

Fix a minimal generating sequence $\x$ for $\m$,
and consider the Koszul complex $K=K^R(\x)$.
The finiteness of $\ciid_R(C)$ and $\pd_R(K)$ implies
$\ciid_R(K\lotimes_R C)<\infty$ by~\cite[(4.4.b)]{sather:cidfc}.
The complex $K\lotimes_R C$ has finite length homology, so~\cite[(3.6)]{sather:cidfc}
implies that there exists a quasi-deformation $R\xra{\vf} R'\xla{\rho} Q$ such that
$Q$ is complete and $\id_Q(R'\lotimes_R(K\lotimes_R C))$ is finite.
Again invoking Lemma~\ref{lem0201}, there is 
a semidualizing $Q$-complex $N$ such that
$R'\lotimes_Q N\simeq R'\lotimes_R C$.  
Lemma~\ref{lem0202} provides 
a sequence $\y\in Q$ such that 
the Koszul complex $L=K^Q(\y)$ 
satisfies $R'\lotimes_QL\simeq R'\lotimes_RK$.
By Lemma~\ref{lem0204} there is an isomorphism 
\begin{align*}
R'\lotimes_R(K\lotimes_R C)
&\simeq \rhom_Q(\rhom_Q(R'\lotimes_QL,Q), N)
\end{align*}
and so the finiteness of 
\begin{align*}
\id_Q(R'\lotimes_R(K\lotimes_R C))
&=\id_Q(\rhom_Q(\rhom_Q(R'\lotimes_QL,Q), N))
\end{align*}
implies $\id_Q(N)<\infty$ by Fact~\ref{ext03'}.
Because $N$ is semidualizing for $Q$, this means that $N$ is dualizing for $Q$, and
Lemma~\ref{lem0203}\eqref{lem0203b}
implies that $C$ is dualizing for $R$.
\end{proof}

\begin{disc} \label{disc0401}
Let $X$ be an $R$-complex.
By combining Fact~\ref{cidim02} with
Property~\ref{picard02c} and
Proposition~\ref{cidim03}\eqref{cidim03item1}, 
one concludes that the following conditions are equivalent:
\begin{enumerate}[\quad\rm(i)]
\item \label{disc0401a}
$X$ is tilting;
\item \label{disc0401c}
$X$ is semidualizing and $\pd_R(X)$ is finite;
\item \label{disc0401b}
$X$ is semidualizing and $\cidim_R(X)$ is finite.
\end{enumerate}
\end{disc}

The next result sets the stage for our equivalence relation on $\s(R)$.

\begin{prop} \label{picard05}
Let  $P$ be a tilting $R$-complex, and let $B$, $C$, and $X$ be  
$R$-complexes with $B$ and $C$ semidualizing.
\begin{enumerate}[\quad\rm(a)]
\item \label{picard05item1}
The $R$-complex $X$ is tilting
(respectively, semidualizing or dualizing)
if and only if $P\lotimes_R X$ is so.
\item \label{picard05item2}
There are isomorphisms
\[ \rhom_R(C,P\lotimes_R C)\simeq P\qquad\text{and}\qquad
\rhom_R(P\lotimes_R C,C)\simeq \rhom_R(P,R). \]
\item \label{picard05item3}
If 
$P\simeq\rhom_R(B,C)$, 
then $C\simeq P\lotimes_R B$.
\end{enumerate}
\end{prop}

\begin{proof}
\eqref{picard05item1}
For each maximal ideal $\m\subset R$, we have
$$(X\lotimes_R P)_{\m}
\simeq X_{\m}\lotimes_{R_{\m}}P_{\m}
\sim X_{\m}\lotimes_{R_{\m}}R_{\m}
\simeq X_{\m}.$$
Hence, Proposition~\ref{picard09} shows that $X$ is tilting
if and only if $P\lotimes_R X$ is tilting;
and $X$ is semidualizing 
if and only if $P\lotimes_R X$ is semidualizing by~\eqref{sp02}.
If $X$ is dualizing, then the previous sentence shows that $P\lotimes_R X$
is semidualizing,
and~\cite[(4.5.F)]{avramov:hdouc} implies $\id_R(P\lotimes_R X)<\infty$;
that is, $P\lotimes_R X$ is dualizing.  Conversely, if
$P\lotimes_R X$ is dualizing, then 
the previous sentence implies that
$\rhom_R(P,R)\lotimes_RP\lotimes_R X$ is dualizing as well,
because $\rhom_R(P,R)$ is tilting; hence, the isomorphism
$$\rhom_R(P,R)\lotimes_RP\lotimes_R X
\simeq X$$
implies that $X$ is dualizing.

\eqref{picard05item2}
Fact~\ref{bass02} implies $P\in\catac(R)$, and this explains
the first isomorphism.
The second isomorphism follows from Hom-tensor adjointness along with the 
assumption $\rhom_R(C,C)\simeq R$.

\eqref{picard05item3}
Assume $P\simeq\rhom_R(B,C)$. 
In particular, we conclude that $\rhom_R(B,C)$ is semidualizing by~\eqref{picard02c},
and so Theorem~\ref{intthmA} implies that $C$ is in $\catbb(R)$.
This explains the first  isomorphism in the following sequence 
\[ C\simeq\rhom_R(B,C)\lotimes_R B \simeq P\lotimes_R B \]
and the second isomorphism is by assumption.
\end{proof}

\begin{disc} \label{picard06}
In conjunction with the basic properties of derived tensor products, 
Proposition~\ref{picard05}\eqref{picard05item1} 
yields a well-defined left $\dpic(R)$-action on $\s(R)$
given by
$[P][C]:=[P\lotimes_R C]$.  Note that the commutativity of 
tensor product shows that the oppositely defined right-action 
is equivalent to this one.

This action yields the following relation on $\s(R)$:  $[B]\approx [C]$ if
$[B]=[P][C]$ for some $[P]\in\dpic(R)$, that is, if $[B]$ is in the orbit
of $[C]$ under the $\dpic(R)$-action.
\end{disc}

We conclude this section with several readily-verified properties of this relation.

\begin{props} \label{picard07}
Let $B$ and $C$ be semidualizing $R$-complexes.
\begin{subprops} \label{picard07d}
For each integer $n$, we have $[\shift^n C]=[(\shift^n R)\lotimes_R C]\approx[C]$.
\end{subprops}
\begin{subprops} \label{picard07a}
The relation $\approx$ is an equivalence relation
because $\s(R)$ is partitioned into orbits by the
$\dpic(R)$-action.
\end{subprops}
\begin{subprops} \label{picard07b}
One has $[C]\approx [R]$ if and only if $[C]\in\dpic(R)$
if and only if $\pd_R(C)<\infty$;
see Remark~\ref{disc0401}.
\end{subprops}
\begin{subprops} \label{picard07c}
When $R$ is local, we have
$[B]\approx[C]$ if and only if $B\sim C$ by~\eqref{picard03}.
\end{subprops}
\end{props}

\section{Detecting Equivalences} \label{sec04}

As we note in the introduction, the results of this section are of the 
following type: Assuming $\ext^n_R(B,C)=0$ for all $n\gg 0$
\emph{and a lot more}, we show $[B]\approx [C]$.  
The first result of this type
is the following version of Property~\ref{sdc02item2}
which also extends Proposition~\ref{picard09} and
Property~\ref{picard07c}. 
Note that the symmetry in
conditions~\eqref{picard08a} and~\eqref{picard08b}
implies that
conditions~\eqref{picard08d} and~\eqref{picard08c}
are symmetric  as well.

\begin{prop} \label{picard08}
Let $B$ and $C$ be semidualizing $R$-complexes.
The following conditions are equivalent:
\begin{enumerate}[\quad\rm(i)]
\item \label{picard08a}
$[B]\approx[C]$; 
\item \label{picard08b}
$B$ is $C$-reflexive and
$C$ is $B$-reflexive;
\item \label{picard08d}
$\ext^n_R(B,C)=0$ for all $n\gg 0$,
and $B_{\m}\sim C_{\m}$ for each maximal ideal $\m\subset R$;
\item \label{picard08c}
$\rhom_R(B,C)$ is a tilting $R$-complex.
\end{enumerate}
When these conditions are satisfied, there are isomorphisms
\begin{align*}
C&\simeq\rhom_R(B,C)\lotimes_RB
&&\text{and}
&B&\simeq\rhom_R(C,B)\lotimes_RC.
\end{align*}
\end{prop}

\begin{proof}
\eqref{picard08a}$\implies$\eqref{picard08b}.
Assume $[B]\approx[C]$ and fix a tilting
$R$-complex $P$ such that $B\simeq P\lotimes_R C$.  
This yields the first isomorphism in
each line of the following display
\begin{gather*}
\rhom_R(B,C)\simeq\rhom_R(P\lotimes_R C,C)\simeq \rhom_R(P,R)\\
\rhom_R(C,B)\simeq\rhom_R(C,P\lotimes_R C)\simeq P
\end{gather*}
and the second isomorphism in each line is from
Proposition~\ref{picard05}\eqref{picard05item2}.
In particular, the complexes $\rhom_R(B,C)$
and $\rhom_R(C,B)$ are homologically bounded.

As in the proof of Proposition~\ref{picard05}\eqref{picard05item1},
we have $B_{\m}\sim C_{\m}$ for each maximal ideal
$\m\subseteq R$.
In particular,  $B_{\m}$ is $C_{\m}$-reflexive and
$C_{\m}$ is $B_{\m}$-reflexive for each $\m$. Using the 
conclusion of the previous paragraph with~\eqref{sdc02item5}
we conclude that $B$ is $C$-reflexive and
$C$ is $B$-reflexive.

\eqref{picard08b}$\implies$\eqref{picard08d}.
Assume that $B$ is $C$-reflexive and
$C$ is $B$-reflexive.  
In particular, this implies $\ext^n_R(B,C)=0$ for all $n\gg 0$.
It follows from~\eqref{sdc02item5}
that $B_{\m}$ is $C_{\m}$-reflexive and
$C_{\m}$ is $B_{\m}$-reflexive for each maximal ideal $\m\subset R$.
Hence, we have $B_{\m}\sim C_{\m}$ for each $\m$
by~\eqref{sdc02item2}.

\eqref{picard08d}$\implies$\eqref{picard08c}.
Assume $\ext^n_R(B,C)=0$ for all $n\gg 0$,
and $B_{\m}\sim C_{\m}$ for each maximal ideal $\m\subset R$.
The first of these assumptions implies that $\rhom_R(B,C)$ is homologically finite.
The second assumption yields the second isomorphism in the next sequence
and the third isomorphism is from~\eqref{sp02}
$$\rhom_R(B,C)_{\m}
\simeq\rhom_{R_{\m}}(B_{\m},C_{\m})
\sim\rhom_{R_{\m}}(C_{\m},C_{\m})
\simeq R_{\m}.$$
Hence $\rhom_R(B,C)$ is tilting by Proposition~\ref{picard09}.

\eqref{picard08c}$\implies$\eqref{picard08a}.
This follows directly from Proposition~\ref{picard05}\eqref{picard05item3}.

When the conditions \eqref{picard08a}--\eqref{picard08c}
are satisfied, the desired isomorphisms
follow from the symmetry of condition~\eqref{picard08b}
using Proposition~\ref{picard05}\eqref{picard05item3}.
\end{proof}

\begin{para} \label{sdc04}
\emph{Proof of Theorem~\ref{intthmB}.}
The equivalence 
\eqref{sdc04item1}$\iff$\eqref{sdc04item6}
is in Proposition~\ref{picard08}, 
while \eqref{sdc04item6}$\implies$\eqref{sdc04item2} is by definition,
and~\eqref{sdc04item2}$\implies$\eqref{sdc04item3}
follows from Fact~\ref{cidim02}.

\eqref{sdc04item3}$\implies$\eqref{sdc04item1}.
Assume that $\cidim_{R}(\rhom_R(B, C))$ is finite.
In particular, we have $\ext_R^n(B,C)=0$ for all $n\gg 0$, and
$$\cidim_{R_{\m}}(\rhom_{R_{\m}}(B_{\m}, C_{\m}))
=\cidim_{R_{\m}}(\rhom_R(B, C)_{\m})<\infty$$
for each maximal ideal $\m\subset R$.  Proposition~\ref{picard08} shows that
it suffices to prove $B_{\m}\sim C_{\m}$ for each $\m$, so we may pass to the
ring $R_{\m}$ and assume that $R$ is local.

As $\cidim_{R}(\rhom_R(B, C))<\infty$
we know from~\cite[(5.1.c)]{sather:cidfc} that
$\rhom_R(B,C)$ is $C$-reflexive.
Using Theorem~\ref{intthmA}, we conclude that 
$\rhom_R(B,C)$ is semidualizing and $B$ is $C$-reflexive.  
Proposition~\ref{cidim03}\eqref{cidim03item1}
and Property~\ref{sdc01item0}
imply
$\rhom_R(B, C) \sim R$, and this 
yields the second isomorphism in  the next sequence
$$B\simeq \rhom_R (\rhom_R (B, C), C) \sim 
\rhom_R(R,C)\simeq C.$$  
The first isomorphism comes from the fact that $B$ is $C$-reflexive.

\eqref{sdc04item1}$\implies$\eqref{sdc04item4}. 
Assume $[B]\approx [C]$, and fix a tilting $R$-complex $P$ 
and an isomorphism
$\alpha\colon B\res P\lotimes_R C$.  
We verify the containment $\catbc(R)\subseteq\catbb(R)$.
Once we do this, 
the symmetry of $\approx$ implies 
$\catbb(R)\subseteq\catbc(R)$,
thus showing $\catbc(R)=\catbb(R)$.

Let $X\in\catbc(R)$.  
The complex $\rhom_R(C,X)$ is homologically bounded,
and so Lemma~\ref{tilt01} implies $\rhom_R(C,X)\in\catb_P(R)$.
Thus, the following sequence shows that $\rhom_R(B,X)$ is also
homologically bounded
$$\rhom_R(B,X)\simeq\rhom_R(P\lotimes_RC,X)
\simeq\rhom_R(P,\rhom_R(C,X))\in\catdb(R).$$
To finish showing that $X$ is in $\catbb(R)$,
we use the next commutative diagram
to conclude that $\xi^B_X$ is an isomorphism.
$$
\xymatrix{
C\lotimes_RP\lotimes_R\rhom_R(P,\rhom_R(C,X)) 
\ar[d]^{\simeq}_{\beta\lotimes_R\sigma} \ar[rr]_-{\simeq}^-{C\lotimes_R\xi^P_{\rhom_R(C,X)}}
&& C\lotimes_R\rhom_R(C,X) \ar[d]_{\simeq}^{\xi^C_X} \\
P\lotimes_RC\lotimes_R\rhom_R(P\lotimes_RC,X) 
\ar[d]^{\simeq}_{P\lotimes_RC\lotimes_R\rhom_R(\alpha,X)} 
&& X \\
P\lotimes_RC\lotimes_R\rhom_R(B,X) \ar[rr]^-{\alpha^{-1}\lotimes_R\rhom_R(B,X)}_-{\simeq}
&& B\lotimes_R\rhom_R(B,X) \ar[u]_-{\xi^B_X}
}
$$
Here, the morphism $\beta$ is tensor-commutativity, and
$\sigma$ is Hom-tensor adjointness.
The morphisms  $\xi^C_X$ and $\xi^P_{\rhom_R(C,X)}$ are isomorphisms
because $X\in\catbc(R)$ and $\rhom_R(C,X)\in\catb_P(R)$.

\eqref{sdc04item4}$\implies$\eqref{sdc04item5}.
This follows from the conditions $B\in\catbb(R)$ and $C\in\catbc(R)$.

\eqref{sdc04item5}$\implies$\eqref{sdc04item1}.
Assume that $B\in\catbc(R)$ and $C\in\catbb(R)$.
Theorem~\ref{intthmA}
implies that $C$ is $B$-reflexive and
$B$ is $C$-reflexive, and so Proposition~\ref{picard08} implies $[B]\approx [C]$.
\qed
\end{para}

Our next result is a version of Theorem~\ref{intthmB} for tensor products.

\begin{prop} \label{sdc05}
Let  $B$ and $C$ be
semidualizing $R$-complexes.
The 
following conditions are equivalent:
\begin{enumerate}[\quad\rm(i)]
\item  \label{sdc05item1}
$B$ and $C$ are both tilting $R$-complexes, i.e, $[B]\approx [R]\approx[C]$;
\item  \label{sdc05item4}
$B\lotimes_R C$ is a tilting $R$-complex;
\item  \label{sdc05item2}
$\pd_R(B\lotimes_R C)<\infty$;
\item  \label{sdc05item3}
$\cidim_R(B\lotimes_R C)<\infty$.
\end{enumerate}
\end{prop}

\begin{proof} 
The implication \eqref{sdc05item1}$\implies$\eqref{sdc05item4}
is in Proposition~\ref{picard05}\eqref{picard05item1}.
Also, \eqref{sdc05item4}$\implies$\eqref{sdc05item2}
is by definition, and
\eqref{sdc05item2}$\implies$\eqref{sdc05item3}
follows from Fact~\ref{cidim02}.

\eqref{sdc05item3}$\implies$\eqref{sdc05item1}.
Assume $\cidim_R(B\lotimes_R C)<\infty$.
For each maximal  $\m\subset R$, we have
$$\cidim_{R_{\m}}(B_{\m}\lotimes_{R_{\m}} C_{\m})
=\cidim_{R_{\m}}((B\lotimes_{R} C)_{\m})<\infty.$$
We show that
$B_{\m}\sim R_{\m}\sim C_{\m}$ for each $\m$,
and so Proposition~\ref{picard09} implies
that $B$ and $C$ are tilting.  Thus, we may assume without loss
of generality that $R$ is local.
Hence, 
there is a quasi-deformation
$R\to R'\from Q$ such that $\pd_Q(R'\lotimes_R (B\lotimes_R C))$ is finite
and $Q$ is complete; see~\cite[Thm.\ F]{sather:cidfc}.  
Because $Q$ is complete, 
Lemma~\ref{lem0201} implies that there  are semidualizing $Q$-complexes
$M$ and $N$ such that
$R'\lotimes_Q M\simeq R'\lotimes_R B$
and $R'\lotimes_Q N\simeq R'\lotimes_R C$.

Each of the following isomorphisms is either standard or by assumption
\begin{align*}
R'\lotimes_R(B\lotimes_R C)
&\simeq(R'\lotimes_R B)\lotimes_{R'}(R'\lotimes_R C) \\
&\simeq(R'\lotimes_Q M)\lotimes_{R'}(R'\lotimes_Q N) \\
&\simeq R'\lotimes_Q (M\lotimes_Q N).
\end{align*}
and so the finiteness of 
$\pd_Q(R'\lotimes_Q (M\lotimes_Q N))
=\pd_Q(R'\lotimes_R(B\lotimes_R C))$
implies $\pd_Q(M)<\infty$ and $\pd_Q(N)<\infty$
by Fact~\ref{ext03'}.
Because $M$ and $N$ are semidualizing for $Q$, this means 
$M\sim Q\sim N$ because of~\eqref{sdc01item0}, and the conclusion
$B\sim R\sim C$ follows
from Lemma~\ref{lem0203}\eqref{lem0203a}.
\end{proof}

Comparing Theorem~\ref{intthmB} and Proposition~\ref{sdc05},
one might feel that there is a missing condition in~\ref{sdc05},
namely $\catab(R)=\catac(R)$.  The next result shows that this 
condition is in fact not  equivalent to those in~\ref{sdc05}.

\begin{prop} \label{sdc05'}
For 
semidualizing $R$-complexes $B$ and $C$,
one has $[B]\approx [C]$ if and only if
$\catab(R)=\catac(R)$.
\end{prop}

\begin{proof}
If $[B]\approx [C]$, then argue
as in the proof of the implication
\eqref{sdc04item1}$\implies$\eqref{sdc04item4} in Theorem~\ref{intthmB} 
to conclude $\catab(R)=\catac(R)$.

For the converse, assume $\catab(R)=\catac(R)$, and fix a 
faithfully injective $R$-module $E$.\footnote{Recall 
that $E$ is \emph{faithfully injective} if 
$\Hom_R(-,E)$ is faithfully exact, that is, if,
for every sequence $S$ of $R$-modules $S$ is exact
if and only if $\Hom_R(S,E)$ is exact.
For example, the $R$-module
$E=\oplus_{\m}E_R(R/\m)$ is faithfully injective,
where the sum is taken over the set of maximal ideals $\m\subset R$;
see~\cite[(3.2.2)]{ishikawa:fef}.}
Observe that there is an isomorphism $\rhom_R(E,E)\simeq \hom_R(E,E)$
because $E$ is injective.  Furthermore, the $R$-module $F=\hom_R(E,E)$
is faithfully flat by~\cite[(1.5)]{ishikawa:oimafm}.
In particular, an $R$-complex $X$ is in $\catdb(R)$  if and only if
$X\lotimes_R F$ is in $\catdb(R)$, and a morphism of $R$-complexes $\alpha$
is an isomorphism in $\D(R)$ if and only if 
$\alpha\lotimes_R F$ is an isomorphism in $\D(R)$.

We show that the containment $\catab(R)\supseteq\catac(R)$
implies $C\in\catbb(R)$.
The reverse containment
then yields $B\in\catbc(R)$,
and so Theorem~\ref{intthmB}  implies $[B]\approx[C]$.

We know that $C$ is in $\catbc(R)$, and hence~\cite[(2.1.f)]{christensen:apac}
implies $\rhom_R(C,E)\in\catac(R)\subseteq\catab(R)$.
It follows from~\cite[(2.1.e)]{christensen:apac} that 
$\rhom_R(\rhom_R(C,E),E)$ is in $\catbb(R)$.
Hom-evaluation~\eqref{basics03} yields the first isomorphism in the
next sequence
$$\rhom_R(\rhom_R(C,E),E)\simeq C\lotimes_R\rhom_R(E,E)\simeq C\lotimes_RF$$
and so $C\lotimes_R F\in\catbb(R)$.
From this, we know that the complex 
$$\rhom_R(B,C\lotimes_R F)\simeq \rhom_R(B,C)\lotimes_R F$$
is homologically bounded; the isomorphism is tensor-evaluation~\eqref{basics03}.
As we noted above, this implies that
$\rhom_R(B,C)$ is homologically bounded.
From the following commutative diagram
$$\xymatrix{
B\lotimes_R \rhom_R(B,C\lotimes_R F) 
\ar[d]_-{\omega_{BCF}}^-{\simeq} \ar[rrd]^-{\xi^B_{C\lotimes_R F}}_{\simeq} \\
B\lotimes_R\rhom_R(B,C)\lotimes_R F
\ar[rr]_-{\xi^B_C\lotimes_R F} 
&& C\lotimes_R F
}$$
we conclude that $\xi^B_C\lotimes_R F$ is an isomorphism.
As we noted above, this implies that $\xi^B_C$ is an isomorphism,
and so $C\in\catbb(R)$, as desired.
\end{proof}

The next  two results  are injective versions
of Theorem~\ref{intthmB} and Proposition~\ref{sdc05}.
Note that we do  not assume that $R$ has a dualizing complex
in either result.

\begin{prop} \label{sdc06}
Let  $B$ and $C$ be
semidualizing $R$-complexes.
The 
following conditions are equivalent:
\begin{enumerate}[\quad\rm(i)]
\item  \label{sdc06item1}
$B$ is a tilting $R$-complex and $C$ is a dualizing $R$-complex;
\item  \label{sdc06item4}
$\rhom_R(B, C)$ is a dualizing $R$-complex;
\item  \label{sdc06item2}
$\id_R(\rhom_R(B, C))<\infty$;
\item  \label{sdc06item3}
$\ciid_R(\rhom_R(B, C))<\infty$.
\end{enumerate}
\end{prop}

\begin{proof} 
We have \eqref{sdc06item4}$\implies$\eqref{sdc06item2}
by definition, and 
\eqref{sdc06item2}$\implies$\eqref{sdc06item3}
is in Lemma~\ref{solve01}\eqref{solve01b}.

\eqref{sdc06item1}$\implies$\eqref{sdc06item4}.
If $B$ is tilting and $C$ is dualizing,
then $\pd_R(B)<\infty$ and $\id_R(C)<\infty$,
and so~\cite[(4.1.I)]{avramov:hdouc} implies $\id(\rhom_R(B,C))<\infty$.
Also  $\rhom_R(B, C)$ is semidualizing by~\eqref{sdc02item3},
and hence is dualizing.

\eqref{sdc06item3}$\implies$\eqref{sdc06item1}.
Assume that $\ciid_R(\rhom_R(B, C))$ is finite.

We first show that $C$ is dualizing for $R$ and $B\sim R$ when $R$ is local.
Fix a minimal generating sequence $\x$ for $\m$,
and consider the Koszul complex $K=K^R(\x)$.
As $\ciid_R(\rhom_R(B, C))$ and $\pd_R(K)$ are finite, \cite[(4.4.b)]{sather:cidfc}
implies
$\ciid_R(K\lotimes_R \rhom_R(B, C))<\infty$.
The complex $K\lotimes_R \rhom_R(B, C)$ has finite length homology, so~\cite[(3.6)]{sather:cidfc}
provides a quasi-deformation $R\xra{\vf} R'\xla{\rho} Q$ such that
$Q$ is complete and $\id_Q(R'\lotimes_R(K\lotimes_R \rhom_R(B, C)))$ is finite.
By Lemmas~\ref{lem0201} and~\ref{lem0202}, 
there are semidualizing $Q$-complexes $M$ and $N$ such that
$R'\lotimes_Q M\simeq R'\lotimes_R B$ and
$R'\lotimes_Q N\simeq R'\lotimes_R C$, and there is
a sequence $\y\in Q$ such that 
the Koszul complex $L=K^Q(\y)$ satisfies
$R'\lotimes_QL\simeq R'\lotimes_RK$.

In the following sequence, the first and third isomorphisms
are combinations of tensor-evaluation~\eqref{basics03}
and Hom-tensor adjointness, and the second isomorphism
is by assumption
\begin{align*}
R'\lotimes_R\rhom_R(B,C)
&\simeq\rhom_{R'}(R'\lotimes_R B,R'\lotimes_R C) \\
&\simeq\rhom_{R'}(R'\lotimes_Q M,R'\lotimes_Q N) \\
&\simeq R'\lotimes_Q \rhom_Q(M,N).
\end{align*}
Hence, Lemma~\ref{lem0204}
provides the following isomorphism
\begin{align*}
R'\lotimes_R(K\lotimes_R \rhom_R(B, C)) 
&\simeq \rhom_Q(\rhom_Q(R'\lotimes_QL,Q), \rhom_Q(M,N)).
\end{align*}
By Fact~\ref{ext03'}, the finiteness of 
\begin{align*}
\id_Q(R'\lotimes_R(K\lotimes_R \rhom_R(B, C))) \hspace{-2cm} \\
&=\id_Q(\rhom_Q(\rhom_Q(R'\lotimes_QL,Q), \rhom_Q(M,N)))
\end{align*}
implies $\id_Q(\rhom_Q(M,N))<\infty$ and 
further $\pd_Q(M)<\infty$ and $\id_Q(N)<\infty$.
Because $M$ and $N$ are semidualizing for $Q$, this 
implies that $N$ is dualizing for $Q$ and
$M\sim Q$. Finally, Lemma~\ref{lem0203}
shows that $C$ is dualizing for $R$ and $B\sim R$.
This concludes the proof when $R$ is local.

We now show that $B$ is a tilting $R$-complex and $C$ is a dualizing $R$-complex
in general.
Our assumptions  guarantee that $\rhom_R(B, C)$ is homologically finite,
and $\ciid_{R_{\m}}(\rhom_{R_{\m}}(B_{\m}, C_{\m}))<\infty$ 
for each maximal ideal
$\m\subset R$.
The local implies $B_{\m}\sim R_{\m}$
and that $C_{\m}$ is dualizing for $R_{\m}$ for each $\m$.
From this, Proposition~\ref{picard09} shows that $B$ is tilting,
and so Lemma~\ref{tilt01} implies $C\in\catbb(R)$.
From Theorem~\ref{intthmA} we conclude that
$\rhom_R(B,C)$ is semidualizing
and so
$\rhom_R(B,C)$ is dualizing
by Proposition~\ref{cidim03}\eqref{cidim03item2}.
In the next sequence, the first isomorphism is 
Hom-tensor adjointness, the second isomorphism is from~\eqref{picard02a},
and the third one is standard
\begin{align*}
\rhom_R(\rhom_R(B,R),\rhom_R(B,C))
&\simeq\rhom_R(\rhom_R(B,R)\lotimes_RB,C) \\
&\simeq\rhom_R(R,C) \\
&\simeq C.
\end{align*}
Since $\pd_R(\rhom_R(B,R))<\infty$ and $\id_R(\rhom_R(B,C))<\infty$,
the first complex in the above sequence has finite injective dimension
by~\cite[(4.1.I)]{avramov:hdouc}. Hence $\id_R(C)<\infty$ and $C$ is dualizing, as desired.
Thus, we localize at $\m$ 
in order to assume that
$R$ is local with maximal ideal $\m$.
\end{proof}

\begin{prop} \label{sdc07}
Let $B$ and $C$ be
semidualizing $R$-complexes.
The 
following conditions are equivalent:
\begin{enumerate}[\quad\rm(i)]
\item  \label{sdc07item1}
$R$ admits a dualizing complex $D$ such that $B\simeq \rhom_R(C,D)$;
\item  \label{sdc07item4}
$B\lotimes_R C$ is a dualizing $R$-complex;
\item  \label{sdc07item2}
$\id_R(B\lotimes_R C)<\infty$;
\item  \label{sdc07item3}
$\ciid_R(B\lotimes_R C)<\infty$.
\end{enumerate}
\end{prop}

\begin{proof} 
We have \eqref{sdc07item4}$\implies$\eqref{sdc07item2}
by definition, and
\eqref{sdc07item2}$\implies$\eqref{sdc07item3}
is in Lemma~\ref{solve01}\eqref{solve01b}.
The implication
\eqref{sdc07item1}$\implies$\eqref{sdc07item4}
follows from~\eqref{sdc02item3}
using standard isomorphisms.

\eqref{sdc07item3}$\implies$\eqref{sdc07item1}.
Assume that $\ciid_R(B\lotimes_R C)$ is finite.
We show that the complex $D=B\lotimes_R C$ is dualizing for $R$;
then~\eqref{sdc02item4}
yields the desired isomorphism:
$$\rhom_R(C,D) =\rhom_R(C,B\lotimes_R C)\simeq B.$$
The finiteness of $\ciid_R(B\lotimes_R C)$ 
implies that $B\lotimes_R C$ is homologically finite,
and $\cidim_{R_{\m}}(B_{\m}\lotimes_{R_{\m}} C_{\m})<\infty$ 
for each maximal ideal
$\m\subset R$.
We  show how this implies that the complex
$(B\lotimes_R C)_{\m}\simeq B_{\m}\lotimes_{R_{\m}}C_{\m}$ 
is dualizing for $R_{\m}$ for each $\m$.
From this, 
Lemma~\ref{solve01}\eqref{solve01c} implies 
$\id_R(B\lotimes_R C)<\infty$,
and~\eqref{sp02} shows that $B\lotimes_R C$
is semidualizing for $R$;
hence $B\lotimes_R C$ is dualizing.
Thus, we localize at $\m$ 
in order to assume that
$R$ is local with maximal ideal $\m$.

The assumption $\ciid_R(B\lotimes_RC)<\infty$ implies
that $B\lotimes_RC$ is homologically finite,
and so $B\lotimes_RC\in\catbc(R)$ by~\cite[(5.2.b)]{sather:cidfc}.
Fact~\ref{bass02} then implies $B\in\catac(R)$, and so
$B\lotimes_R C$ is semidualizing by Corollary~\ref{semidual6}.
It follows from Proposition~\ref{cidim03}\eqref{cidim03item2}
that $B\lotimes_R C$ is dualizing as desired.
\end{proof}

\section{Ext-vanishing, Tor-vanishing and Cohen-Macaulayness} \label{sec02}

The following is proved in~\cite[(1.3)]{avramov:edcrcvct} and serves as
our motivation for this section.

\begin{fact} \label{ABS01}
Let $R$ be a local ring admitting a dualizing complex $D$ and assume $\inf(D)=0$.
If $\ext^n_R(D,R)=0$ for $i=1,\ldots,\dim(R)$, then $R$ is Cohen-Macaulay.
\end{fact}

We  generalize this fact to the realm of semidualizing complexes 
after the following lemma which compliments~\eqref{dual04'}.

\begin{lem} \label{lem0601}
Let $R$ be a local ring and $C$ a semidualizing $R$-complex.
If $\amp(C)=\cmd(R)$ and 
$\supp_R(\HH_{\sup(C)}(C))=\spec(R)$, then $C$ is Cohen-Macaulay.
\end{lem}

\begin{proof}
After replacing $C$ with
$\shift^{-\inf(C)}C$, we assume without loss of generality that $\inf(C)=0$,
which implies $\sup(C)=\amp(C)=\cmd(R)$.
We show that, for each $\p\in\spec(R)$, we have
$\dim(R/\p)-\inf(C_{\p})\leq\depth_R(C)$.
From this it follows that
$$\depth_R(C)\leq\dim_R(C)
=\sup\{\dim(R/\p)-\inf(C_{\p})\mid\p\in\spec(R)\}
\leq\depth_R(C)$$
and so $C$ is Cohen-Macaulay.

Fix a prime $\p\in\spec(R)$.
Our assumptions yield the following sequence
$$0=\inf(C)
\leq \inf(C_{\p})\leq\sup(C_{\p})
=\sup(C)=\cmd(R).$$
The complex $C_{\p}$ is semidualizing for $R_{\p}$
and so~\eqref{dual04'} implies
$$\cmd(R)-\inf(C_{\p})
=\amp(C_{\p})\leq\cmd(R_{\p}).$$
This explains the first inequality in the following sequence
\begin{align*}
\dim(R/\p)-\inf(C_{\p})
&\leq\dim(R/\p)-\cmd(R)+\cmd(R_{\p})\\
&=[\dim(R/\p)+\dim(R_{\p})-\dim(R)]+\depth(R)-\depth(R_{\p})\\
&\leq\depth(R)-\depth(R_{\p})\\
&\leq\depth(R)\\
&=\depth_R(C).
\end{align*}
The first equality is by definition;
the second and third inequalities are standard;
and the final equality is from~\cite[(3.2.a)]{christensen:scatac}.
\end{proof}

See Remark~\ref{ABS02} for an explicit discussion of the connection
between Fact~\ref{ABS01} and the next result. 

\begin{prop} \label{ext01} 
Let $R$ be a local ring and fix  semidualizing $R$-complexes $B$ and $C$
such that $\inf(B)=0=\inf(C)$. 
\begin{enumerate}[\quad\rm(a)]
\item \label{ext01item0}
Fix a prime  $\p\in\spec(R)$
such that $R_{\p}$ is Cohen-Macaulay, e.g., $\p\in\Min(R)$.  
If $i=\sup(B_{\p})$ and $j=\sup(C_{\p})$,
then $\ext^{i-j}_R(B,C)\neq 0$.
\item \label{ext01item2}
Fix an integer $s\geq \sup(C)$.
If $\ext^n_R(B,C)=0$ for $n=-s+1,\ldots,\dim(R)$, then 
$B$ is isomorphic to a module in $\D(R)$, $s= \sup(C)$ and 
$\supp_R(\HH_{s}(C))=\spec(R)$.
\item \label{ext01item2'}
If $\ext^n_R(B,C)=0$ for $n=-\cmd(R)+1,\ldots,\dim(R)$, then 
$C$ is Cohen-Macaulay.
\item \label{ext01item3}
Assume that $R$ admits a dualizing complex $D$ such that $\inf(D)=0$
and that $\sup(C) = 0$. If 
$\ext^n_R(\rhom_R(C,D),C)=0$ for $n=1,\ldots,\dim(R)$, then 
$R$ is Cohen-Macaulay.
\item \label{ext01item4}
Assume that $R$ admits a dualizing complex $D$ such that $\inf(D)=0$
and that $B$ is Cohen-Macaulay. If
$\ext^n_R(B,\rhom_R(B,D))=0$ for $n=1,\ldots,\dim(R)$, then 
$R$ is Cohen-Macaulay.
\end{enumerate}
\end{prop}

\begin{proof}
\eqref{ext01item0}
The equality $\supp_R(B)=\spec(R)$ from~\eqref{sp01}
implies $0\leq i\leq\sup(B)$, and similarly $0\leq j\leq\sup(C)$.
Furthermore, since 
$R_{\p}$ is Cohen-Macaulay, we have $\amp(B_{\p})=0=\amp(C_{\p})$
by~\eqref{sdc01item2}.
In particular, there are isomorphisms
$$B_{\p}\simeq\shift^i\HH_i(B)_{\p} \qquad \qquad \qquad
C_{\p}\simeq\shift^j\HH_j(C)_{\p} $$
and these provide the second isomorphism in the next sequence
\begin{align*}
\ext^{i-j}_R(B,C)_{\p}
&\cong \ext^{i-j}_{R_{\p}}(B_{\p},C_{\p})
\cong \ext^{i-j}_{R_{\p}}(\shift^i\HH_i(B)_{\p},\shift^j\HH_j(C)_{\p}) \\
&\cong \ext^0_{R_{\p}}(\HH_i(B)_{\p},\HH_j(C)_{\p}) 
\cong \Hom_{R_{\p}}(\HH_i(B)_{\p},\HH_j(C)_{\p}) 
\neq 0.
\end{align*}
The first, third and fourth isomorphisms are standard;
and the nonvanishing holds 
by~\cite[(3.6)]{holm:fear}
because $\HH_i(B)_{\p}$ 
is a semidualizing $R_{\p}$-module.
We conclude that $\ext^{i-j}_R(B,C)\neq 0$, as desired.

\eqref{ext01item2}
Fix a prime ideal $\p\in\Min(R)$, and set
$i=\sup(B_{\p})$ and $j=\sup(C_{\p})$.
The assumption $\inf(B)=0=\inf(C)$ implies
$0\leq i\leq\sup(B)$ and $0\leq j\leq\sup(C)\leq s$.
This justifies the first five inequalities in the next sequence
$$-s\leq-\sup(C)\leq-j\leq i-j\leq i\leq\sup(B)=\amp(B)\leq\cmd(R)\leq\dim(R).$$
The sixth inequality is in~\eqref{dual04'}, and
the equality is by assumption.

From part~\eqref{ext01item0} 
we know $\ext^{i-j}_R(B,C)\neq 0$. Hence, the previously displayed
inequalities along with our vanishing hypothesis imply 
$$-s=-\sup(C)=-j= i-j$$
Hence, we have $j=\sup(C)=s$ and $i=0$. The first equality 
with the definition $j=\sup(C_{\p})$ implies $\p\in\supp_R(\HH_{s}(C))$.
Since $\p$ is an arbitrary prime in $\Min(R)$, it follows that
$\Min(R)\subseteq\supp_R(\HH_{s}(C))$.
Because $\supp_R(\HH_{s}(C))$ is Zariski closed in $\spec(R)$,
the equality $\supp_R(\HH_{s}(C))=\spec(R)$ follows. 

As in the previous paragraph, the equality $0=i=\sup(B_{\p})$ for each
$\p\in\Min(R)$ implies $\supp_R(\HH_{0}(B))=\spec(R)$.  
In particular, for each $\q\in\spec(R)$, we have $\inf(B_{\q})=0$.
For each
$\q\in\ass_R(\HH_{\sup(B)}(B))$
this yields the second equality in the next sequence
$$\sup(B)=\inf(B_{\q})=0.$$
The first equality is from~\eqref{sp02}.
Thus $B$ is
isomorphic to a module in $\D(R)$. 

\eqref{ext01item2'}
Since $\cmd(R)\geq \amp(C)=\sup(C)$ by~\eqref{dual04'},
part~\eqref{ext01item2}
implies $\sup(C)=\cmd(R)$ and $\supp_R(\HH_{\sup(C)}(C))=\spec(R)$.
Now apply Lemma~\ref{lem0601} to conclude that $C$ is
Cohen-Macaulay.

\eqref{ext01item3}
From~\eqref{sdc02item3} we know $\inf(\rhom_R(C,D))=0$, and so
part~\eqref{ext01item2} provides the first equality in the next sequence
$$0=\amp((\rhom_R(C,D))=\cmd(C)=\cmd(R).$$
The remaining equalities are from~\eqref{sdc02item3} and~\eqref{dual04'}, 
respectively, using the assumption 
$\amp(C)=0$.
Hence, $R$ is Cohen-Macaulay, as desired.

\eqref{ext01item4}
If $B$ is Cohen-Macaulay, then $\amp(\rhom_R(B,D))=0$ 
by~\eqref{sdc02item3}.
The isomorphism $\rhom_R(\rhom_R(B,D),D)\simeq B$ shows that we may apply 
part~\eqref{ext01item3} using the complex $C=\rhom_R(B,D)$ to conclude
that $R$ is Cohen-Macaulay.
\end{proof}

The next remark shows the need for the  hypotheses 
in Proposition~\ref{ext01}\eqref{ext01item3}
and~\eqref{ext01item4}.

\begin{disc} \label{ABS02'}
Let $R$ be a local ring and fix  semidualizing $R$-complexes $B$ and $C$
such that $\inf(B)=0=\inf(C)$ and
$\ext^n_R(B,C)=0$ for $n=1,\ldots,\dim(R)$.
In view of Fact~\ref{ABS01} and Proposition~\ref{ext01},
one might be tempted to conclude that 
$R$ is Cohen-Macaulay.
However, the example $B=R=C$ shows that this conclusion need not 
follow.  Using the same $B$ and $C$ with
Proposition~\ref{tor01} in mind, we see that 
the assumption 
$\tor_n^R(B,C)=0$ for $n=1,\ldots,2\dim(R)$ also does not necessarily
imply that $R$ is Cohen-Macaulay.
\end{disc}

We next explain why Fact~\ref{ABS01}
is a special case of Proposition~\ref{ext01}.

\begin{disc} \label{ABS02}
Let $R$ be a local ring, and let $D$ be a dualizing complex 
for $R$ such that $\inf(D)=0$ and
$\ext^n_R(D,R)=0$ for $i=1,\ldots,\dim(R)$.
Notice that $D$ is Cohen-Macaulay.
Apply Proposition~\ref{ext01}\eqref{ext01item2} with $B=D$ and
$C=R$ to conclude
$\sup(D)=0$.  Hence, the equality $\sup(D)=\cmd(R)$ implies that
$R$ is Cohen-Macaulay.  One can also apply
Proposition~\ref{ext01}\eqref{ext01item3}
with $C=R$ or Proposition~\ref{ext01}\eqref{ext01item4} with $B=D$
to draw the same conclusion.
\end{disc}

The proof of the following is very similar to that of Proposition~\ref{ext01}.

\begin{prop} \label{tor01}
Let $R$ be a local ring and fix  semidualizing $R$-complexes $B$ and $C$
such that $\inf(B)=0=\inf(C)$. Set $r=\sup(B)$ and $s=\sup(C)$.
\begin{enumerate}[\quad\rm(a)]
\item \label{tor01item1}
If $\p\in\ass_R(\HH_{r}(B))$
and  $j=\inf(C_{\p})$
then $\tor_{r+j}^R(B,C)\neq 0$.
\item \label{tor01item2}
If $\tor_n^R(B,C)=0$ for $n=1,\ldots,2\dim(R)$, then 
$B$ and $C$ are isomorphic to modules in $\D(R)$.
\item \label{tor01item3}
Assume that either $B$ or $C$  is Cohen-Macaulay.
If $\tor_n^R(B,C)=0$ for $n=1,\ldots,2\dim(R)$, then 
$R$ is Cohen-Macaulay.
\qed
\end{enumerate}
\end{prop}

\section{Special Cases of Question~\ref{tach03}\eqref{tach03item1}} \label{sec06}

We open this section with an example showing why we need 
$R$ to be local in Questions~\ref{tach03}\eqref{tach03item1}
and~\ref{tpq02}\eqref{tpq02item1}.

\begin{ex} \label{local01}
Let $(R_1,\m_1)$ and $(R_2,\m_2)$ be commutative noetherian local rings, and set
$R=R_1\times R_2$.  Assume that $R_2$ is Cohen-Macaulay and non-Gorenstein
and admits a dualizing module $D_2$.
Fix an integer $m>\dim(R)+1$, and consider
the $R$-complex $B=R_1\oplus \shift^mD_2$.

Set $\m'=\m_1\times R_2$ and
$\m''=R_1\times \m_2$.  Recall that $\{\m',\m''\}$ is precisely the
set of maximal ideals of $R$.
Also, there are isomorphisms 
\begin{align*}
R_{\m'}&\cong(R_1)_{\m_1}\cong R_1
&R_{\m''}&\cong(R_2)_{\m_2}\cong R_2\\
B_{\m'}&\simeq(R_1)_{\m_1}\simeq R_1
&B_{\m''}&\simeq\shift^m(D_2)_{\m_2}\simeq\shift^mD_2.
\end{align*}
These isomorphisms with~\eqref{sp02} imply that 
$B$ is semidualizing for
$R$.  

We claim that $\ext^n_R(B,R)=0$ for $n=1,\ldots,m-1$.
Indeed, the above isomorphisms imply
\begin{align*}
\ext^n_R(B,R)_{\m'}
&\cong\ext^n_{R_{\m'}}(B_{\m'},R_{\m'})
\cong\ext^n_{R_1}(R_1,R_1)=0
&&\text{for $n\geq 1$} \\
\ext^n_R(B,R)_{\m''}
&\cong\ext^n_{R_{\m''}}(B_{\m''},R_{\m''})\\
&\cong\ext^n_{R_2}(\shift^mD_2,R_2)
\cong\ext^{n-m}_{R_2}(D_2,R_2)=0
&&\text{for $n<m$.} 
\end{align*}
This establishes the claim since it can be verified locally.

Suppose that $B$ is $R$-reflexive.  
We conclude from~\eqref{sdc02item5}
that $B_{\m''}$ is $R_{\m''}$-reflexive, that is
$D_2$ is $R_2$-reflexive.  From~\eqref{sdc02item3} we know that
$R_2$ is $D_2$-reflexive, so~\eqref{sdc02item2} implies $D_2\sim R_2$,
contradicting the assumption that $R_2$ is non-Gorenstein;
see~\eqref{dual03item2}.
This shows why we must assume that 
$R$ is local in Question~\ref{tach03}\eqref{tach03item1}.

For Question~\ref{tpq02}\eqref{tpq02item1} argue as above
to conclude that 
$\tor^n_R(B,B)=0$ for $n=1,\ldots,2m-1$.
Suppose that $B\lotimes_R B$ is semidualizing for $R$. 
From~\eqref{sp02} we conclude that
that $(B\lotimes_R B)_{\m''}$ is semidualizing for $R_{\m''}$, that is
$D_2\lotimes_{R_2}D_2$ is semidualizing for $R_2$.  Using~\cite[(3.2)]{frankild:sdcms}
we conclude $D_2\sim R_2$, again
contradicting the assumption that $R_2$ is non-Gorenstein.
\end{ex}

\begin{para} \label{ABS2.1}
\emph{Proof of Theorem~\ref{intthmC}.}
In each case, Proposition~\ref{ext01} implies that $R$ is Cohen-Macaulay.
In particular~\eqref{dual04'} implies 
$$\sup(B)=\sup(C)=\sup(D)=0=\amp(\rhom_R(B,D))=\amp(\rhom_R(C,D)).$$  
Hence, we assume that $B$, $C$ and $D$
are modules.
From~\eqref{sdc02item3} there are equalities
$\inf(\rhom_R(B,D))=\inf(\rhom_R(C,D))=0$
and this yields isomorphisms 
$$\rhom_R(B,D)\simeq\Hom_R(B,D)
\qquad\text{and}
\qquad\rhom_R(C,D)\simeq\Hom_R(C,D).$$

\eqref{ABS2.1item1}
We claim that there is an isomorphism
\begin{equation} \label{homeval01} \tag{\ref{ABS2.1}.1}
\Hom_R(\Hom_R(B,\Hom_R(B,D)),D)\cong B\otimes_R B.
\end{equation}
The argument is akin to that of~\cite[(B.3)]{avramov:edcrcvct}.
Let $D\res J$ be a minimal injective resolution, and set 
$G=\rhom_R(B,\Hom_R(B,D))$.
In particular, we have $J_p=0$ when $p\geq 1$ and when $p<-\dim(R)$.
From~\cite[(B.1)]{avramov:edcrcvct}, there is a 
strongly convergent spectral 
sequence
$$E^2_{p,q}=\HH_p(\Hom_R(\HH_{-q}(G),J))\implies 
\HH_{p+q}(\Hom_R(G,J)).$$
The differentials act in the pattern $d^r_{p,q}\colon E^r_{p,q}\to E^r_{p-r,q+r-1}$.

Our vanishing assumptions for $\ext^q_R(B,\Hom_R(B,D))\cong\HH_{-q}(G)$
and $J_p$ imply $E^2_{p,q}=0$ when $1\leq q\leq\dim(R)$,
when $p\geq 1$ and when $p<-\dim(R)$.  
In particular, we have $E^{\infty}_{-q,q}=0$ when $q\neq 0$
and 
$$E^{\infty}_{0,0}=E^2_{0,0}=\HH_0(\Hom_R(\HH_{0}(G),J))
\cong\Hom_R(\Hom_R(B,\Hom_R(B,D)),D).$$
The isomorphism is by the left-exactness of 
$\Hom_R(\Hom_R(B,\Hom_R(B,D)),-)$.
The vanishing $E^{\infty}_{-q,q}=0$ when $q\neq 0$ yields the first isomorphism 
in the following
\begin{equation} \label{2.1isos}  \tag{\ref{ABS2.1}.2}
\HH_0(\Hom_R(G,J))\cong E^{\infty}_{0,0}
\cong\Hom_R(\Hom_R(B,\Hom_R(B,D)),D).
\end{equation}
On the other hand, the first isomorphism in the next sequence
is by construction
\begin{align*}
\Hom_R(G,J)
&\simeq \rhom_R(\rhom_R(B,\Hom_R(B,D)),D)\\
&\simeq B\lotimes_R\rhom_R(\rhom_R(B,D),D)\\
&\simeq B\lotimes_R B.
\end{align*}
The second isomorphism is by Hom-evaluation~\eqref{basics03},
and the third is from~\eqref{sdc02item3}.
Taking degree-0 homology  in this sequence 
and using Fact~\ref{ext03}, we have
$$\HH_0(\Hom_R(G,J))\cong
\HH_0(B\lotimes_R B)\cong B\otimes_R B.$$
With~\eqref{2.1isos}, this provides the isomorphism~\eqref{homeval01}.

The isomorphism~\eqref{homeval01} yields the first equality
in the following sequence
\begin{align*}
\ass_R(B\otimes_R B)
&=\ass_R(\Hom_R(\Hom_R(B,\Hom_R(B,D)),D))\\
&=\supp_R(\Hom_R(B,\Hom_R(B,D)))\cap\ass_R(D)\\
&\subseteq\ass_R(R).
\end{align*}
The second equality is standard, and the containment is 
by~\eqref{dual04''}.
Now use~\cite[(2.2)]{avramov:edcrcvct} to conclude $B\cong R$.

\eqref{ABS2.1item0}
Because $\amp(C)=0$, we know that $\rhom_R(C,D)$ 
is Cohen-Macaulay by~\eqref{sdc02item3}.
Since $\rhom_R(\rhom_R(C,D),D)\simeq C$  we may apply 
part~\eqref{ABS2.1item1} with $B=\rhom_R(C,D)$ 
to conclude
$\rhom_R(C,D)\simeq R$ and hence
\begin{xxalignat}{3}
  &{\hphantom{\square}}& C
  &\simeq\rhom_R(\rhom_R(C,D),D)\simeq\rhom_R(R,D)\simeq D. &&\qed
\end{xxalignat}
\end{para}

\begin{disc} \label{disc-noTor}
It is natural to try to use
Theorem~\ref{intthmC} 
to answer Question~\ref{tpq02} when $R$ is generically Gorenstein. 
To see what goes wrong,
let $R$, $B$, $C$, and $D$ be as in Theorem~\ref{intthmC},
and assume that 
$\tor_n^R(B,C)=0$ for $n=1,\ldots, d$ for some integer $d$.
Our Tor-vanishing hypothesis combines with
the isomorphism
$$\rhom_R(B\lotimes_R C,D)\simeq\rhom_R(B,\rhom_R(C,D))$$
to imply
$\ext^n_R(B,\rhom_R(C,D))=0$ for 
$d-\dim(R)$ consecutive values of $n\geq 1$.
However, we cannot apply
Theorem~\ref{intthmC} 
to $B$ and $\rhom_R(C,D)$ even when $d$ is large,
because the Ext-vanishing
begins at $n=\dim(R)+1$, not at $n=1$.
\end{disc}

To highlight the difficulties outlined in the previous remark, we pose the following
question which is a special case of Question~\ref{tpq02}\eqref{tpq02item2}.

\begin{question} \label{moretpq}
Assume that $R$ is a complete
Cohen-Macaulay local ring of dimension $d\geq 1$ that is Gorenstein
on the punctured spectrum.
Let $B$ and $C$ be semidualzing $R$-modules
such that $\tor_n^R(B,C)=0$ for all $n\geq 1$.
Must $B\lotimes_R C$ be semidualizing for $R$?
\end{question}

In the next result, the special
case $B=D$  in part~\eqref{ABS3.1item1} 
(or $C=R$ in part~\eqref{ABS3.1item2})
is exactly~\cite[(3.1)]{avramov:edcrcvct}.
Note that, if $Q$ admits a dualizing complex, that is, if 
$Q$ is a homomorphic image of a Gorenstein ring, or if $Q$ is excellent, then $Q$ 
automatically has
Gorenstein formal fibres.

\begin{prop} \label{ABS3.1}
Let $Q$ be a generically Gorenstein, local ring with Gorenstein
formal fibres.  Let $\mathbf{x}=x_1,\ldots,x_c\in Q$ be a $Q$-regular sequence,
and set $R=Q/(\mathbf{x})$.
Assume that $R$  admits a dualizing complex $D$ such that $\inf(D)=0$.
Let $B$ and $C$ be  semidualizing $R$-complexes such that $\inf(B)=0=\inf(C)$.
Assume that $B$ is Cohen-Macaulay and $\sup(C)=0$.
\begin{enumerate}[\quad\rm(a)]
\item \label{ABS3.1item1}
If  $\ext^n_R(B,\rhom_R(B,D))=0$
for $n=1,\ldots,\dim(R)+1$, then $B\simeq R$.
\item \label{ABS3.1item2}
If  $\ext^n_R(\rhom_R(C,D),C)=0$
for $n=1,\ldots,\dim(R)+1$, then $C\simeq D$.
\end{enumerate}
\end{prop}

\begin{proof}
In each case, Proposition~\ref{ext01} implies that $R$ is Cohen-Macaulay
and hence, so is $Q$.
As in the proof of Theorem~\ref{intthmC} 
we assume that $B$, $C$ and $D$
are modules and we have
$\rhom_R(B,D)\simeq\Hom_R(B,D)$
and $\rhom_R(C,D)\simeq\Hom_R(C,D)$.

\eqref{ABS3.1item1}
We first reduce to the case where $c=1$, as in the proof of~\cite[(3.1)]{avramov:edcrcvct}.
From~\cite[(3.2)]{avramov:edcrcvct} there is a $Q$-regular sequence
$\mathbf{y}\in Q$ such that $(\mathbf{y})Q=(\mathbf{x})Q$ and
$P=Q/(y_1,\ldots,y_{c-1})$ is generically Gorenstein.  Note that $P$ has 
Gorenstein formal fibres because $Q$ does, and so we may replace $Q$ with $P$
in order to assume $c=1$.

Next, we reduce to the case where
$Q$ (and hence $R$) is complete.
The sequence $\mathbf{x}=x_1$ is $\comp{Q}$-regular because
$\comp{Q}$ is $Q$-flat, and
there is an isomorphism $\comp{R}\cong \comp{Q}/(x_1)$.
The ring $\comp{R}$ is Cohen-Macaulay with dualizing 
module $\comp{R}\otimes_R D$ by~\eqref{sp03}.  
Furthermore, the module
$\comp{R}\otimes_R B\simeq\comp{R}\lotimes_R B$ is $\comp{Q}$-semidualizing 
by~\eqref{sp03}.
The standard base-change result for flat extensions yields the following
isomorphisms
\begin{align*}
\ext^n_{\comp{R}}(\comp{R}\otimes_R B,\Hom_{\comp{R}}(\comp{R}\otimes_R B,\comp{R}\otimes_R D))
&\cong\ext^n_{\comp{R}}(\comp{R}\otimes_R B,\comp{R}\otimes_R \Hom_R(B,D)) \\
&\cong\comp{R}\otimes_R \ext^n_R(B,\Hom_R(B,D))
\end{align*}
for each integer $n$. In particular, we have
$\ext^n_{\comp{R}}(\comp{R}\otimes_R B,\Hom_{\comp{R}}(\comp{R}\otimes_R B,\comp{R}\otimes_R D))
=0$ for $i=1,\ldots,\dim(R)+1$, that is, $i=1,\ldots,\dim(\comp{R})+1$.
If it follows that $\comp{R}\otimes_R B\cong\comp{R}$, then~\eqref{sp03}
implies $B\cong R$, as desired.

To finish the reduction, we need to explain why $\comp{Q}$ is generically Gorenstein.
Because $\comp{Q}$ is Cohen-Macaulay, there is an equality
$\ass(\comp{Q})=\Min(\comp{Q})$.  Fix a prime ideal $\q\in\Min(\comp{Q})$ and
set $\p=\q\cap Q$. The going-down-property implies $\p\in\Min(Q)$, and so
$Q_{\p}$ is Gorenstein by assumption. The induced morphism
$Q_{\p}\to \comp{Q}_{\q}$ is flat and local, and 
the fact that $Q$ has Gorenstein formal fibres implies that this map
has Gorenstein closed fibre.  It follows that $\comp{Q}_{\q}$ 
is Gorenstein, as desired.

Now, we assume that $Q$ and $R$ are complete,
so the ring $Q$ admits
a dualizing module $D^Q$.  
Lemma~\ref{lem0201}
implies that $Q$
admits a semidualizing complex $A$ such that $B\simeq R\lotimes_Q A$,
and~\eqref{sp03} implies $\amp(A)=\amp(B)=0$.
That is, $A$ is a semidualizing $Q$-module such that
$B\cong R\otimes_Q A$.
(This is also proved by Gerko~\cite[(3)]{gerko:osmagi}.)
Note that  $x_1$ is $A$-regular and $D^Q$-regular
by~\eqref{dual04''}, and also $D\cong R\otimes_Q D^Q$.  
By~\cite[(3.3.1)]{avramov:edcrcvct} the hypothesis
$\ext^n_R(B,\Hom_R(B,D))=0$
for $i=1,\ldots,\dim(R)+1$ implies
$\ext^n_P(A,\Hom_R(A,D^Q))=0$
for $i=1,\ldots,\dim(Q)$, and so Theorem~\ref{intthmC} implies
$A\cong Q$.  It follows that
$B\cong R\otimes_Q A\cong R\otimes_Q Q\cong R$,
as desired.

\eqref{ABS3.1item2}
This follows from part~\eqref{ABS3.1item1} as in the proof
of Theorem~\ref{intthmC}\eqref{ABS2.1item0}.
\end{proof}

Our final result follows from Proposition~\ref{ABS3.1} using the argument
of~\cite[(4.1), (4.2)]{avramov:edcrcvct} which is, in turn,
the special
case $B=D$  in part~\eqref{ABS4.1item1} 
(or $C=R$ in part~\eqref{ABS4.1item2}).

\begin{prop} \label{ABS4.1}
Let $R$ be a  local ring of the form $Q/\fa$,
where $Q$ is a Gorenstein local ring and $\fa$ is an ideal for which there is a 
sequence of links $\fa\sim\fb_1\sim\cdots\sim\fb_s\sim\fb$ with $\fb$ a generically
complete intersection ideal.  (This is the case, e.g., if $R$ is a homomorphic image
of a Gorenstein ring and $\operatorname{codepth}(R)\leq 2$.)
Let $D$ be  a dualizing $R$-complex such that $\inf(D)=0$,
and let $B$ and $C$ be  semidualizing $R$-complexes such that $\inf(B)=0=\inf(C)$.
Assume that $B$ is Cohen-Macaulay and $\sup(C)=0$.
\begin{enumerate}[\quad\rm(a)]
\item \label{ABS4.1item1}
If  $\ext^n_R(B,\rhom_R(B,D))=0$
for $n=1,\ldots,\dim(R)+1$, then $B\simeq R$.
\item \label{ABS4.1item2}
If  $\ext^n_R(\rhom_R(C,D),C)=0$
for $n=1,\ldots,\dim(R)+1$, then $C\simeq D$.
\qed
\end{enumerate}
\end{prop}

\appendix
\section{Complete Intersection Dimensions and Quasi-deformations}
\label{sec07}

This appendix contains several technical lemmas
for use in Sections~\ref{sec05} and~\ref{sec04}.

\begin{lem} \label{solve01}
Let $X$ be a homologically finite $R$-complex.
\begin{enumerate}[\quad\rm(a)]
\item \label{solve01b}
There is an inequality $\ciid_R(X)\leq\id_R(X)$ with equality when $\id_R(X)<\infty$.
\item \label{solve01c}
If $\ciid_R(X)<\infty$ and $\id_{R_{\m}}(X_{\m})<\infty$ for each maximal ideal
$\m\subset R$, then $\id_R(X)=\ciid_R(X)<\infty$.
\end{enumerate}
\end{lem}

\begin{proof}
Assume that $\id_{R_{\m}}(X_{\m})<\infty$ for every maximal ideal 
$\m\subset R$.
Fact~\ref{cidim02}
shows $\ciid_{R_{\m}}(X_{\m})<\infty$ for each $\m$,
and so the following Bass formulas 
come from~\cite[(2.11)]{sather:cidfc} and~\cite[(2.7.b)]{foxby:cim}
$$\ciid_{R_{\m}}(X_{\m})=\depth_{R_{\m}}(X_{\m})
-\inf(X_{\m})=\id_{R_{\m}}(X_{\m}).$$
This yields the second equality in the next sequence
\begin{align*}
\id_R(X)
&=\sup\left\{\id_{R_{\m}}(X_{\m}) \mid \text{$\m\subset R$ is a maximal ideal}\right\} \\
&=\sup\left\{\ciid_{R_{\m}}(X_{\m}) \mid \text{$\m\subset R$ is a maximal ideal}\right\} \\
&=\ciid_R(X).
\end{align*}
The first 
equality is from~\cite[(5.3.I)]{avramov:hdouc},
and the third equality is in Definition~\ref{cidim01'}.
Part~\eqref{solve01c} follows immediately.

\eqref{solve01b} Assume that $\id_R(X)<\infty$.
This implies $\id_{R_{\m}}(X_{\m})<\infty$ for every maximal 
$\m\subset R$
by~\cite[(5.3.I)]{avramov:hdouc}, so the desired conclusion
follows from the previous display.
\end{proof}

\begin{lem} \label{lem0201}
Assume that $R$ is local and let $R\to R'\from Q$ be a quasi-deformation
such that $Q$ is complete. For each semidualizing $R$-complex $C$
there exists a semidualizing $Q$-complex $N$ such that
$R'\lotimes_QN\simeq R'\lotimes_R C$.
\end{lem}

\begin{proof}
The complex $R'\lotimes_R C$ is semidualizing for $R'$
by~\eqref{sp03}. 
In particular, this implies $\ext^2_{R'}(R'\lotimes_RC,R'\lotimes_RC)=0$.
Thus, because $Q$ is complete and the map $Q\to R'$ is 
surjective with kernel generated by a $Q$-sequence, 
a lifting result of Yoshino~\cite[(3.2)]{yoshino} provides a homologically
finite $Q$-complex $N$ such that $R'\lotimes_QN\simeq R'\lotimes_R C$.
Again using the fact that $R'\lotimes_R C$ is semidualizing for $R'$,
this isomorphism implies that $N$ is semidualizing for $Q$ by~\cite[(4.5)]{frankild:rrhffd}.
\end{proof}

In the next result, the Koszul complex over $R$ on a sequence $\x$
is denoted $K^R(\x)$.

\begin{lem} \label{lem0202}
Assume that $R$ is local and let $R\xra{\vf} R'\xla{\rho} Q$ 
be a quasi-deformation.
For each sequence $\x=x_1,\ldots,x_n\in R$ there exists a sequence
$\y=y_1,\ldots,y_n\in Q$ such that 
$R'\lotimes_RK^R(\x)\simeq R'\lotimes_QK^Q(\y)$.
\end{lem}

\begin{proof}
Because $\rho$ is 
surjective, there exist elements $y_i\in Q$ such that $\rho(y_i)=\vf(x_i)$
for each $i$.  The desired isomorphism now follows from the next sequence
\begin{xxalignat}{3}
  &{\hphantom{\square}}& R'\lotimes_RK^R(\x)
  &\simeq K^{R'}(\phi(\x))\simeq K^{R'}(\rho(\y)) 
     \simeq R'\lotimes_QK^Q(\y). &&\qedhere
\end{xxalignat}
\end{proof}

\begin{lem} \label{lem0203}
Assume that $R$ is local and let $R\to R'\from Q$ 
be a quasi-deformation.
Let $C$ be a semidualzing $R$-complex and assume that
$N$ is a semidualzing $Q$-complex such that 
$R'\lotimes_QN\simeq R'\lotimes_R C$.
\begin{enumerate}[\quad\rm(a)]
\item \label{lem0203a}
If $N\sim Q$, then $C\sim R$.
\item \label{lem0203b}
If $N$ is dualizing for $Q$, then $C$ is dualizing for $R$.
\end{enumerate}
\end{lem}

\begin{proof}
\eqref{lem0203a}
Assuming $N\sim Q$, the next isomorphisms follow from our hypotheses
$$R'\lotimes_R C\simeq R'\lotimes_Q N\sim R'\lotimes_Q Q\simeq R'\simeq R'\lotimes_R R$$
and so~\eqref{sp03} implies $C\sim R$.

\eqref{lem0203b}
Assume that $N$ is dualizing for $Q$.
Because the map 
$Q\to R'$ is surjective with kernel generated by a $Q$-regular sequence,
we conclude from~\eqref{sp03} 
that $R'\lotimes_Q N\simeq R'\lotimes_R C$ is dualizing for $R'$, and 
similarly that $C$ is dualizing for $R$.
\end{proof}

\begin{lem} \label{lem0204}
Assume that $R$ is local and let $R\to R'\from Q$ 
be a quasi-deformation.
Let $K$ and $X$ be $R$-complexes and let $L$ and $Y$
be $Q$-complexes such that 
$R'\lotimes_Q L\simeq R'\lotimes_R K$ and 
$R'\lotimes_Q Y\simeq R'\lotimes_R X$.
If $L$ is homologically finite over $Q$ and $\pd_Q(L)$ is finite,
then
$$R'\lotimes_R(K\lotimes_R C)
\simeq \rhom_Q(\rhom_Q(R'\lotimes_QL,Q), N).$$
\end{lem}

\begin{proof}
Because $L$ and $R'$ are both homologically finite $Q$-complexes
of finite projective dimension,
Fact~\ref{ext03'} implies
$\pd_Q(R'\lotimes_QL)<\infty$.
Also, the fact that 
$R'\lotimes_QL$ is homologically finite over $Q$ implies 
$\pd_Q(\rhom_Q(R'\lotimes_QL,Q))<\infty$ by~\cite[(2.13)]{christensen:scatac}.

Our assumptions yield the second isomorphism in the next sequence
while the first and third isomorphisms are standard.
\begin{align*}
R'\lotimes_R(K\lotimes_R X)
&\simeq (R'\lotimes_RK)\lotimes_{R'}(R'\lotimes_R X) \\
&\simeq (R'\lotimes_QL)\lotimes_{R'}(R'\lotimes_QY) \\
&\simeq (R'\lotimes_QL)\lotimes_Q Y \\
&\simeq \rhom_Q(\rhom_Q(R'\lotimes_QL,Q),Q)\lotimes_Q Y \\
&\simeq \rhom_Q(\rhom_Q(R'\lotimes_QL,Q), Y)
\end{align*}
The fourth isomorphism is from~\eqref{sdc02item1}
and the fifth one is in~\cite[(1.7.b)]{frankild:rrhffd};
these rely on the finiteness of
$\pd_Q(R'\lotimes_QL)$ and $\pd_Q(\rhom_Q(R'\lotimes_QL,Q))$.
\end{proof}

\section*{Acknowledgments}

We are grateful to Lars W.\ Christensen
and Hans-Bj\o rn Foxby
for stimulating
conversations about this research.
We thank the referee for
thoughtful suggestions.


\providecommand{\bysame}{\leavevmode\hbox to3em{\hrulefill}\thinspace}
\providecommand{\MR}{\relax\ifhmode\unskip\space\fi MR }
\providecommand{\MRhref}[2]{%
  \href{http://www.ams.org/mathscinet-getitem?mr=#1}{#2}
}
\providecommand{\href}[2]{#2}

\end{document}